\documentclass[12pt,reqno,a4paper]{amsproc}
\usepackage{amssymb}
\usepackage{amsfonts}
\usepackage{amsmath}
\usepackage{graphicx}
\usepackage{graphics}
\usepackage{color}
\newtheorem{thm}{Theorem}[section]
\newtheorem{lm}{Lemma}[section]
\newtheorem{defn}{Definition}[section]
\newtheorem{pro}{Proposition}[section]
\newtheorem{re}{Remark}[section]
\newtheorem{cor}{Corollary}[section]

\newtheorem{exa}{Example}[section]
\newtheorem{con}{Conjecture}[section]
\usepackage{esint}

\DeclareMathOperator{\tr}{tr}

\DeclareMathOperator{\ric}{Ric}

\DeclareMathOperator{\vol}{vol}
\DeclareMathOperator{\diam}{diam}
\DeclareMathOperator{\Rm}{Rm}

\makeatletter
\newsavebox\myboxA
\newsavebox\myboxB
\newlength\mylenA

\newcommand*\xoverline[2][0.75]{%
    \sbox{\myboxA}{$\m@th#2$}%
    \setbox\myboxB\null
    \ht\myboxB=\ht\myboxA%
    \dp\myboxB=\dp\myboxA%
    \wd\myboxB=#1\wd\myboxA
    \sbox\myboxB{$\m@th\overline{\copy\myboxB}$}
    \setlength\mylenA{\the\wd\myboxA}
    \addtolength\mylenA{-\the\wd\myboxB}%
    \ifdim\wd\myboxB<\wd\myboxA%
       \rlap{\hskip 0.5\mylenA\usebox\myboxB}{\usebox\myboxA}%
    \else
        \hskip -0.5\mylenA\rlap{\usebox\myboxA}{\hskip 0.5\mylenA\usebox\myboxB}%
    \fi}
\makeatother

\begin{document}

\title[]{On the Fill-in of Nonnegative Scalar Curvature Metrics}

\author{Yuguang Shi}
\address [Yuguang Shi] {Key Laboratory of Pure and Applied Mathematics, School of Mathematical Sciences, Peking University, Beijing, 100871, P.\ R.\ China}
\email{ygshi@math.pku.edu.cn}
\thanks{Yuguang Shi, Guodong Wei and Jiantian Zhu are partially supported by NSFC 11671015 and 11731001.  Wenlong Wang is partially supported by NSFC 11671015, 11701326 and BX201700007.}

\author{Wenlong Wang}
\address [Wenlong Wang] {Department of Mathematics and LPMC, Nankai University, Tianjin, 300071, P. R. China}
 \email{wangwl@nankai.edu.cn}

\author{Guodong Wei}
\address [Guodong Wei]  {Key Laboratory of Pure and Applied Mathematics, School of Mathematical Sciences, Peking University, Beijing, 100871, P.~R.~China}
\email{weiguodong@amss.ac.cn}

\author{Jintian Zhu}
\address [Jintian Zhu]   {Key Laboratory of Pure and Applied Mathematics, School of Mathematical Sciences, Peking University, Beijing, 100871, P.~R.~China}
\email{zhujt@pku.edu.cn}

\renewcommand{\subjclassname}{
 \textup{2010} Mathematics Subject Classification}

\date{June 4, 2019}

\begin{abstract}
In the first part of this paper, we consider the problem of  fill-in of nonnegative scalar curvature (NNSC) metrics for a triple of Bartnik data $(\Sigma,\gamma,H)$. We prove that given a metric $\gamma$ on $\mathbf{S}^{n-1}$ ($3\leq n\leq 7$), $(\mathbf{S}^{n-1},\gamma,H)$ admits no fill-in of NNSC metrics provided the prescribed mean curvature $H$ is large enough (Theorem \ref{Thm: no fillin nonnegative scalar 2}). Moreover, we prove that if $\gamma$ is a positive scalar curvature (PSC) metric isotopic to the standard metric on $\mathbf{S}^{n-1}$, then the much weaker condition that the total mean curvature $\int_{\mathbf S^{n-1}}H\,\mathrm d\mu_\gamma$ is large enough rules out NNSC fill-ins, giving an partially affirmative answer to a conjecture by Gromov (see P.\,23 in \cite{Gromov4}).  In the second part of this paper, we investigate the $\theta$-invariant of Bartnik data and obtain some sufficient conditions for the existence of PSC fill-ins.
\end{abstract}
\maketitle
\markboth{Shi Yuguang, Wang Wenlong, Wei Guodong and Zhu Jintian }{Fill-in of NNSC}

\section{Introduction}

A triple of (generalized) Bartnik data $(\Sigma^{n-1},\gamma,H) $ consists of an orientable closed null-cobordant Riemannian manifold $(\Sigma^{n-1},\gamma)$ and a given smooth function $H$ on $\Sigma^{n-1}$. One  basic problem in Riemannian geometry is to study (see \cite{Gromov2}): {\it under what kind of conditions does the Bartnik data  $(\Sigma^{n-1},\gamma,H) $ admit a fill-in metric $g$ with scalar curvature bounded below by a given constant? That is, there are a compact Riemannian manifold $(\Omega^n , g )$ with boundary of scalar curvature $R_g\geq \sigma>-\infty$,  and an isometry $X: (\Sigma^{n-1}, \gamma)\mapsto (\partial \Omega^n,  g|_{\partial \Omega^n})$ so that $H=H_g \circ X$ on $\Sigma$, where $H_g$ is the mean curvature of $\partial \Omega^n$ in $(\Omega^n, g)$ with respect to the outward unit normal vector. }

Note that the above definition of fill-in is different from that in \cite{J}. In our case, if $(\Omega^n, g,X)$ is a fill-in of $(\Sigma^{n-1},\gamma,H)$, we have $\partial \Omega^n=X(\Sigma^{n-1})$ rather than $X(\Sigma^{n-1})\subset\partial \Omega^n$ and $\partial \Omega^n \setminus X(\Sigma^{n-1})$ is allowed to be a closed (possibly disconnected) minimal hypersurface of $(\Omega^n, g)$.
By the gluing arguments in \cite{Schoen-Yau2} and \cite{GL}, it is easy to see our definition is more restrictive than that in \cite{J}.

On the other hand, in \cite{Shi-Tam1} (also see an improvement in \cite{Shi-Tam3}), the first author and his collaborator proved the positivity of Brown-York mass introduced by Brown and York (\cite{BY1,BY2}).

\begin{thm}\label{postiveBYmass1}
	Let $(\Omega^3, g)$ be a $3$-dimensional compact Riemannian manifold with nonnegative scalar curvature and with strictly mean-convex boundary $\partial \Omega$ that consists of spheres with positive Gauss curvature. Then for each component $\Sigma_\ell \subset \partial \Omega$, $\ell = 1, \ldots, k$, 
\begin{equation*} \label{e-shi-tam}
		\mathfrak{m}_{BY}(\Sigma_\ell ;\,\Omega,g)\ge0.
	\end{equation*}
Here $\mathfrak{m}_{BY}(\Sigma_\ell ;\,\Omega,g)$ is the Brown-York mass of $\Sigma_\ell$ in $(\Omega,g)$ defined by
\begin{equation*} \label{e-shi-tam}
		\mathfrak{m}_{BY}(\Sigma_\ell ;\,\Omega,g)=\int_{\Sigma_\ell}\left(H_0-H\right)\,\mathrm{d}\mu,
	\end{equation*}
where $H_0$ is the mean curvature of $\Sigma_\ell$ when isometrically embedded in $\mathbf{R}^3$, and $H$ is the mean curvature of $\Sigma_\ell$ in $(\Omega,g)$. 
	Moreover, equality holds for some $\ell$ if and only if $\partial  \Omega$ has only one component and $(\Omega, g)$ is isometric to a domain in $\mathbf{R}^3$.
\end{thm}

Later, they got a more general result, namely

\begin{thm}\label{postiveBYmass2}
Let $(\Omega^3, g)$ be a $3$-dimensional compact Riemannian manifold
with smooth boundary that is a topological sphere. Suppose the scalar curvature of $(\Omega,g)$ satisfies $R_g\geq -6\kappa^2$, the Gauss curvature of its boundary $\Sigma$ satisfies $K> -\kappa^2$, and the mean curvature $H$ of $\Sigma$ is positive. Then
$$\int_{\Sigma}\left(H_0
-H\right)\cosh\kappa r\,\mathrm{d}\mu\ge0, 
$$
where $H_0$ is the mean curvature of $\Sigma$ when isometrically embedded in $\mathbf{H}^3_{-\kappa^2}$ and $r$ is a geodesic distance defined in (2.1) in \cite{Shi-Tam2}.
\end{thm}

Theorem \ref {postiveBYmass1}, Theorem \ref {postiveBYmass2}, as well as Miao's work \cite{Miao1} and Mantoulidis-Miao's work \cite{MM} imply that for a sufficiently large function $H$, it is impossible  to fill in $(\Sigma^2, \gamma, H)$ with a metric $g$ of $R_g\geq \sigma$ for some constant $\sigma$. Indeed, positivity of Brown-York mass is a necessary condition for $(\Sigma^2, \gamma, H)$ to admit a fill-in  metric $g$ with NNSC provided the Gauss curvature of $(\Sigma^2, \gamma)$ is positive; but it is not sufficient (see \cite{J, JMT} for details).

However, all of above works are mainly suitable for the three  dimensional case; not so many results are known for higher dimensional manifolds, which are obviously worth studying. In \cite{Gromov4}, Gromov proposed the following conjecture (see P.\,23):

\begin{con}
Let  $(X, g)$ be a compact Riemannian manifold with scalar curvature $R\geq \sigma$. Then
$$
\int_{\partial X}H \leq \Lambda, 
$$
where $H$ is the mean curvature of the boundary $\partial X$ in $(X, g)$ with respect to the outward unit normal vector, and $\Lambda$ is a constant depending only on $\sigma$ and  the intrinsic geometry of $(\partial X,g|_{\partial X})$.
\end{con}

One goal of this paper is to give a partially affirmative answer to the above conjecture; see Theorem \ref{nonexistfillin1} below. Before stating our results, we first introduce some notations and conventions. 

Throughout this paper, when we discuss a smooth manifold $\Sigma$, we always fix a differential structure $\mathcal{U}$ on it. Any metric $\gamma$ will be understood to be given by metric component functions on the coordinate charts in  $\mathcal{U}$. That is, we will distinguish between a metric $\gamma$ and its pull back $\phi^* \gamma$ by a diffeomorphism $\phi$. In particular, we always consider $\mathbf{S}^{n-1}$ as the unit sphere in the Euclidean space $\mathbf{R}^n$ with the induced differential structure. We use $\gamma_{std}$ to denote the standard metric on $\mathbf{S}^{n-1}$ induced from the Euclidean space. For $k\geq 2$ and a closed manifold $\Sigma^{n-1}$, let $\mathcal M^k(\Sigma^{n-1})$ be the space of all $C^k$ metrics on $\Sigma^{n-1}$ with the $C^k$-topology. Similarly, let $\mathcal M^\infty(\Sigma^{n-1})$ be the space of all smooth metrics on $\Sigma^{n-1}$ with the $C^\infty$-topology. We define $\mathcal M^k_{psc}\left(\Sigma^{n-1}\right)=\left\{\gamma\in\mathcal M^k\left(\Sigma^{n-1}\right)\,|\,R_\gamma>0\right\}$, and call two metrics $\gamma_0$, $\gamma_1$ in $\mathcal M^k_{psc}(\Sigma^{n-1})$ {\it isotopic} to each other if there exists a continuous path $\gamma: [0,1] \mapsto \mathcal M^k_{psc}(\Sigma^{n-1})$ such that $\gamma(0)=\gamma_0$ and $\gamma(1)=\gamma_1$. Finally we make a convention that unless otherwise specified, the mean curvature of a boundary component is with respect to the unit outer normal.

Our first main result is
\begin{thm}\label{nonexistfillin1}
For $3\leq n\leq 7$ and $k\geq 5$, let $\gamma$ be a smooth metric isotopic to $\gamma_{std}$ in $\mathcal M^k_{psc}(\mathbf S^{n-1})$. Then there exists a constant $h_0=h_0(\gamma)$ such that $(\mathbf S^{n-1},\gamma,H)$ admits no fill-in of nonnegative scalar curvature whenever
\begin{equation*}
H>0\quad\text{\rm and}\quad \int_{\mathbf S^{n-1}}H\,\mathrm d\mu_\gamma>h_0.
\end{equation*}
\end{thm}

 Due to \cite{Marqurs}, we know that any PSC metric $\gamma$ on $\mathbf{S}^3$ is isotopic to $\gamma_{std}$. By Proposition 2.1 and its proof in \cite{CAM}, we may assume the path is smooth. Hence, Theorem \ref{nonexistfillin1} holds for all PSC metrics on $\mathbf{S}^3$. 
 
 For general $\gamma$, we also investigate the same problem. Define
\begin{equation*}
\mathcal M^n_{c,\,d}:=\left\{\gamma\in\mathcal M^\infty\left(\mathbf{S}^{n-1}\right)\,\bigg |
\begin{array}{c}
|\Rm_\gamma|\leq c,\,\diam(\gamma)\leq d\\
\vol(\gamma)=\vol(\gamma_{std})
\end{array}\right\}.
\end{equation*}
 We have

\begin{thm}\label{Thm: no fillin nonnegative scalar 2}
For $3\leq n\leq 7$, given positive constants $c$ and $d$, there exists a universal constant $H_0=H_0(n,c,d)$ such that $(\mathbf S^{n-1},\gamma,H)$ admits no fill-in of nonnegative scalar curvature for any $\gamma\in\mathcal M^n_{c,\,d}$ and $H>H_0$.
\end{thm}

A similar result, which is a sharp pointwise comparison of the mean curvature of the boundary,  was obtained in \cite{Gromov2} (P.\,3); however, the domain of fill-in is assumed to be spin there. 

The assumption $3\leq n\leq 7 $ in Theorem \ref{nonexistfillin1} and Theorem \ref{Thm: no fillin nonnegative scalar 2} is only due to the positive mass theorem, which was claimed to be true for all dimension now (\cite{Schoen-Yau1}); thus, the above two theorems are true for all dimensions not less than three.

Inspired by \cite{Gromov1} (P.\,53--54), for a triple of Bartnik data $(\Sigma^{n-1},\gamma,H)$,  we consider the set of fill-ins $\mathcal{F}=\{(\Omega^{n},g,X)\}$, and define the following $\theta$-invariant of $(\Sigma^{n-1},\gamma,H)$ by
$$\theta(\Sigma^{n-1},\gamma,H)=\sup_{\mathcal{F}}\inf_{\Omega}R_g.$$

Obviously, $\theta$-invariant has deep relations with above fill-in problem. For instance, if $\theta(\Sigma^{n-1},\gamma,H)>\sigma$, then $(\Sigma^{n-1},\gamma,H)$ admits a fill-in with a metric $g$ of $R_g\geq \sigma$; and  if $\theta(\Sigma^{n-1},\gamma, H)<\sigma$, then $(\Sigma^{n-1},\gamma,h)$ admits no fill-in with a metric $g$ of $R_g\geq \sigma$.  

 By some known results, we do have a few examples for this invariant. 
 Due to Theorem 4 in \cite{HW}, we have
 
  \begin{exa}\label{2dexplicitformula}
Let $0\leq H<1$ be a constant. Then
\begin{equation*}
\theta(\mathbf{S}^1,\gamma_{std},H)=2\left(1-H^2\right),
\end{equation*} and it is achieved only by filling in $(\mathbf{S}^1,\gamma_{std},H)$ with a spherical cap of $\mathbf{S}^2_{\frac{1}{\sqrt{1-H^2}}}$, the round sphere of radius $\frac{1}{\sqrt{1-H^2}}$. 
 \end{exa}

 Theorems \ref{postiveBYmass1} and \ref{postiveBYmass2} jointly imply

 \begin{exa}\label{3dexplicitformula}
Let $H\geq 2$ be a constant. Then
\begin{equation*}
\theta(\mathbf{S}^{2},\gamma_{std},H)=6\left(1-\frac{H^2}{4}\right),
\end{equation*}
and it is achieved only by filling in $(\mathbf{S}^{2},\gamma_{std},H)$ with a geodesic ball of $\mathbf{R}^3$ ($H=2$) or $\mathbf{H}^3_{1-H^2/4}$ ($H>2$). This result can be generalized to high dimensions. 
\end{exa}

In \cite{MS}, Mantoulidis-Schoen proved
\begin{exa}
For any metric $\gamma$ on $\mathbf{S}^2$ with $\lambda_1(-\Delta_\gamma+K_\gamma)> 0$, where $K_\gamma$ is the Gauss curvature
of $\gamma$, $\theta(\mathbf{S}^{2},\gamma,0)>0$. 
\end{exa}

According to the counterexample to the Min-Oo's conjecture constructed in \cite{BMN}, we have:
\begin{exa}
For $n\geq 3$, $$\theta(\mathbf{S}^{n-1},\gamma_{std},0)>n(n-1),$$ 
so it is not achieved by the hemisphere with the standard metric. 
\end{exa}

Besides these examples,  we also know some information about the $\theta$-invariant of several special cases. For example, we know $$\theta\left(\mathbf{S}^{n_1}_{\sqrt{n_1}}\times\mathbf{S}^{n_2}_{\sqrt{n_2}}\times\cdots\times\mathbf{S}^{n_k}_{\sqrt{n_k}},\gamma_{can},0\right)\geq
n.$$ Here $\gamma_{can}$ is the product metric of which the $n_i$-th factor is the round metric of radius $\sqrt{n_i}$, and $n=\sum^k_{i=1}{n_i}+1$. In particular, $$\theta(\mathbf{T}^{n-1},\gamma_{can},0)\geq n.$$

 But the $\theta$-invariant is far from being  studied systematically. So, it is important to investigate some basic properties of this invariant. In the following, we always assume $(\Sigma^{n-1},\gamma)$ is a $(n-1)$-dimensional orientable closed null-cobordant Riemannian manifold. We first consider a fill-in that achieves the invariant, called  an extremal fill-in, and get

\begin{thm}\label{staticity}
If $\theta(\Sigma^{n-1},\gamma,H)\geq 0$, then any extremal fill-in of $(\Sigma^{n-1},\gamma,H)$ (if exists) is static.
\end{thm}

Another important feature of the  $\theta$-invariant is the following  monotonicity formula

\begin{thm}\label{Thm: weak monotonicity sigma invariant}
Let $H_1$ and $H_2$ be two functions on $\Sigma^{n-1}$. If $H_1\geq H_2$, then
$\theta(\Sigma^{n-1},\gamma,H_1)\leq \theta(\Sigma^{n-1},\gamma,H_2)$.
\end{thm}

\begin{re}
It should be interesting to see what happens when $\theta(\Sigma^{n-1},\gamma,H_1)=\theta(\Sigma^{n-1},\gamma,H_2)$ in Theorem \ref{Thm: weak monotonicity sigma invariant}. Unfortunately, we cannot address this problem for the time being. 
\end{re}

Since the $\theta$-invariant is monotonically non-increasing with respect to the prescribed mean curvature, the limit at positive infinite mean curvature exists. For round spheres, the limits are negative infinity, and for general Riemannian manifolds we have:

\begin{thm}\label{Thm: Exponential Decay sigma Invariant}
Let $(\Sigma^{n-1},\gamma)$ be a Riemnanian manifold with $R_{\gamma}\geq 0$. Then either
\begin{enumerate}
\item[(1)] for any constant $H$, $$\theta(\Sigma^{n-1},\gamma,H)=+\infty,$$
\item[(2)] or there exist positive constants $H_0$, $C$ depending only on $(\Sigma^{n-1},\gamma)$ and a dimensional constant $\beta>0$ such that for any constant $H\geq H_0$, 
\begin{equation*}
\theta(\Sigma^{n-1},\gamma,H)\leq CH^{-\beta}.
\end{equation*}
\end{enumerate}
\end{thm}

 We have used the idea of torical symmetrization (see \cite{Gromov3}) in the proof of above theorem. We have the following corollary: 

\begin{cor}\label{limofsigma}
Let $(\Sigma^{n-1},\gamma)$ be a Riemnanian manifold with $R_{\gamma}\geq 0$. Then either
\begin{enumerate}
\item[(1)] for any function $H$, $\theta(\Sigma^{n-1},\gamma,H)=+\infty$,  or
\item[(2)] for any $\sigma>0$, there is a constant $H_0>0$, such that $\theta(\Sigma^{n-1}, \gamma,H)<\sigma$ for all $H\geq H_0$.
\end{enumerate}
\end{cor}

By a  gluing argument, we obtain:
\begin{thm}\label{nonnegative}
For  $(\Sigma^{n-1},\gamma)$ with $R_{\gamma}>0$, either
\begin{enumerate}
\item[(1)] $\theta(\Sigma^{n-1},\gamma,0)\geq\min R_\gamma $, or
\item[(2)] $\theta(\Sigma^{n-1},\gamma,0)=0$ and it can not be attained.
\end{enumerate}
\end{thm}

Due to the  arguments in \cite{MM2}, we may construct the so called ``Schwarzschild neck" of $(\Sigma^{n-1},\gamma)$ (see Definition \ref{Schwarzschild-neck}). And by gluing such a neck to a suitable fill-in region, we can obtain some estimates for $\theta(\Sigma^{n-1},\gamma,H)$ with small positive $H$. Namely,
\begin{thm}\label{fillin1}
Let $(\Sigma^{n-1},\gamma,H)$ be a triple of Bartnik data with $H\geq 0$ and $
R_\gamma>\frac{n-2}{n-1}\max H^2$. Then one of the following two alternatives holds:
\begin{enumerate}
\item[(1)] $\theta(\Sigma^{n-1},\gamma,H)\geq\min R_\gamma-\frac{n-2}{n-1}\max H^{2}$.
\item[(2)] $\theta(\Sigma^{n-1},\gamma,H)=0$ and it can not be attained.
\end{enumerate}
\end{thm}

As mentioned above, we are interested in when $\theta(\Sigma^{n-1},\gamma,H)>0$. One possible way is to consider connected components of the set of PSC metrics on $\Sigma^{n-1}$.

\begin{thm}\label{Thm: main1}
Let $\gamma_0$ and $\gamma_1$ be two smooth metrics in $\mathcal
M^k _{psc}(\Sigma^{n-1})$ isotopic to each other. If $(\Sigma^{n-1},\gamma_1,0)$ admits a fill-in of positive scalar curvature,
then for any function $H$ with
\begin{equation}\label{Eq: main1}
H<\left(\frac{n-1}{n-2}\min R_{\gamma_0}\right)^{\frac{1}{2}},
\end{equation}
we have $\theta (\Sigma^{n-1},\gamma_0,H)>0$, namely $(\Sigma^{n-1},\gamma_0,H)$ admits a fill-in of positive scalar curvature. 
\end{thm}

Combining Theorem \ref{Thm: main1} with the result in \cite{Marqurs}, we see that for any metric $\gamma\in\mathcal M^k _{psc}(\mathbf{S}^3)$ and $H$ satisfying \eqref{Eq: main1}, $\theta (\mathbf{S}^3,\gamma,H)>0$.

The main idea to prove Theorem \ref{nonexistfillin1}   and   Theorem \ref{Thm: no fillin nonnegative scalar 2} is to construct an asymptotically flat (see Definition \ref{AFmanifold}) end with NNSC and with $(\mathbf S^{n-1},\gamma,H)$ being an inner boundary; then we show that the ADM mass will be negative provided $\int_{\mathbf S^{n-1}}H\,\mathrm d\mu_\gamma$ or $H$ is large enough (see \eqref{mass}); finally by the positive mass theorem for AF manifolds with corners (see Theorem 1 in \cite{Miao1}), we see that such Bartnik data admits no fill-in  of  NNSC metrics. In order to prove Theorem \ref{staticity}, we first observe that if an extremal fill-in is not static, then due to Theorem 1 in \cite{Cor}, we may raise the scalar curvature of an interior subregion but preserves the metric near the boundary by a compact perturbation of the metric. To get the contradiction, we then have to raise the scalar curvature near the boundary but keep the induced metric on the boundary. We achieve this by doing twice suitable conformal deformations and using Theorem 5 in \cite{BMN}. Via a similar approach, we prove Theorem \ref{Thm: weak monotonicity sigma invariant}. By rescaling and gluing a ``nearly extremal'' fill-in region to a certain neck, we find that the $\theta$-invariant decreases for a certain portion after we raise the mean curvature (see Proposition \ref{Prop: strict inequality 2}); then by an iteration argument, we get the  proof of Theorem \ref{Thm: Exponential Decay sigma Invariant}.

The rest of the paper run as follows: in Section 2 we present some useful lemmas and propositions; in Section 3 we prove the main theorems.

\section{Some Basic Lemmas}
Let us begin with the following notion. 
\begin{defn}\label{AFmanifold}

Let $n\geq 3$. A Riemannian manifold $(M^n, g)$ is said to be asymptotically flat (AF) if there is a compact set $K\subset M^n$ such that $M^n\setminus K$ is diffeomorphic to the exterior of a ball in $\mathbf{R}^n$ and in this coordinate $g$ satisfies
$$
|g_{i j}-\delta_{i j}|+|x|\left|\partial g_{i j}\right|+|x|^{2}\left|\partial^{2} g_{i j}\right|+|x|^{3}\left|\partial^{3} g_{i j}\right|
=O\left(|x|^{-p}\right)$$
for some $p>\frac{n-2}2$. Furthermore, we require that
$$\int_{M^n}|R_g|\,\mathrm{d}\mu_g<\infty.$$
The  Arnowitt-Deser-Misner (ADM) mass \cite{ADM} of $(M^n,g)$ is defined by
\begin{equation*}
m_{ADM}(M^n,g)=\lim_{r\to\infty}\frac{1}{{ 2(n-1)}\omega_{n-1}}\int_{S_r}
\left(g_{ij,i}-g_{ii,j}\right)\nu^j\,\mathrm dS_r,
\end{equation*}
where $S_{r}$ is the coordinate sphere near the infinity, $\nu$ is the Euclidean outward unit normal to $S_r$, and $dS_r$ is the Euclidean area element on $S_r$.
\end{defn}

In the sequel, we are going to construct an AF end with a continuous path in $\mathcal M^k_{psc}(\mathbf S^{n-1})$ with endpoint $\gamma_{std}$. Let $\gamma_0 \in \mathcal M^k_{psc}(\mathbf S^{n-1})$ and $\{\gamma(t)\}_{t\in [0,1]}$ be a continuous path in $\mathcal M^k_{psc}(\mathbf S^{n-1})$ with $\gamma(0)=\gamma_0$, $\gamma(1)=\gamma_{std}$. Without loss of generality, we may assume  $\gamma(t)\equiv\gamma_{std}$ for $t \in [\frac56, 1 ]$. By Proposition 2.1 and its proof in \cite{CAM}, we may also assume the path is smooth. We first have the following lemma:

\begin{lm}\label{Lem: base space for Shi-Tam flow}
For $k\geq 2$, let $\gamma_t:[0,1]\to \mathcal M^k(\mathbf S^{n-1})$ be a smooth path with $\gamma_0=\gamma$ and $\gamma_t\equiv \gamma_{std}$ for $t \in [\frac56, 1]$. Given any $\epsilon>0$, there exists a positive constant $s_0=s_0(\epsilon,\sup_{t\in [0,1]}\|\gamma_t'\|_{\gamma_t})$ such that we can find a $C^k$ metric $\bar g$ on $\mathbf S^{n-1}\times[1,+\infty)$ that has the form
\begin{equation*} 
\bar g=\mathrm ds^2+s^2\tilde\gamma_s,
\end{equation*}
where $\tilde\gamma_s:[1,\infty)\to \mathcal M^k(\mathbf S^{n-1})$ is a smooth path with $\tilde\gamma_1=\gamma$, $\tilde \gamma_{s}\equiv \gamma_{std}$ for $s\geq s_0$,
and satisfies
\begin{equation}\label{Eq: mean curvature estimates on base space}
\left\|\bar A_s-\frac{1}{s}\bar \gamma_s\right\|_{\bar\gamma_s}\leq \frac{\epsilon}{s},
\end{equation}
where $\bar\gamma_s=s^2\tilde \gamma_s$ and $\bar A_s$ is the second fundamental form of the slice $\Sigma_s:=\mathbf S^{n-1}\times\{s\}$ with respect to $\bar g$ and the $\partial_s$-direction. Moreover, the scalar curvatures $R_{\bar\gamma_s}$ and $R_{\bar g}$ are bounded by universal constants depending only on $\epsilon$, $\|R_{\gamma_t}\|_{L^\infty([0,1])}$, $\sup_{t\in [0,1]}\|\gamma_{t}'\|_{\gamma_t}$ and $\sup_{t\in [0,1]}\|\gamma_{t}''\|_{\gamma_t}$.
\end{lm}

\begin{re}
Let $E$ denote $\mathbf{S}^{n-1}\times [1, \infty)$. Obviously, $(E, \bar g)$ is AF; indeed, it is Euclidean for $s>s_0$. 
\end{re}

\begin{proof}[Proof of Lemma \ref{Lem: base space for Shi-Tam flow}]
With $\delta>0$ to be determined later, we define $t:[1,+\infty)\to [0,1)$ by
\begin{equation*}\label{Eq: function t(s)}
t(s)=\frac{2}{\pi}\arctan\left(\delta\ln s\right)
\end{equation*}
and let $\tilde\gamma_{s}=\gamma_{t(s)}$. It is clear that $\tilde\gamma_1=\gamma$ and $\tilde \gamma_s\equiv \gamma_{std}$ for $s\geq s_0$ with
\begin{equation}\label{Eq: expression for s_0}
s_0=\exp\left(\frac{1}{\delta}\tan\big(\frac{5\pi}{12}\big)\right).
\end{equation}
Let $\bar g=\mathrm ds^2+s^2\tilde \gamma_s$. Then $\bar g$ is a $C^k$ metric on $\mathbf S^{n-1}\times [1,+\infty)$.  It is not hard to see
\begin{equation*}
\left\|\bar A_s-\frac{1}{s}\bar \gamma_s\right\|_{\bar\gamma_s}=\frac{\delta}{\pi s(1+\delta^2\ln^2s)}\|\gamma'_{t(s)}\|_{\gamma_{t(s)}}.
\end{equation*}
Therefore, we can choose $\delta$ small enough, depending only on $\epsilon$ and $\sup_{t\in [0,1]}\|\gamma_t'\|_{\gamma_t}$, to obtain (\ref{Eq: mean curvature estimates on base space}). And it follows from (\ref{Eq: expression for s_0}) that $s_0$ depends only on $\epsilon$ and $\sup_{t\in [0,1]}\|\gamma_t'\|_{\gamma_t}$ as well. Finally, the bounds on $R_{\bar\gamma_s}$ and $R_{\bar g}$ come from a straightforward calculation.
\end{proof}

\begin{exa}
Let $\Sigma_0$ be a smooth closed strictly convex hypersurface
in $\mathbf{R}^n$ and $r$ be the distance function to $\Sigma_0$. Then the metric on the exterior region of $\Sigma_0$ is given by $dr^2 + g_r$, where $g_r$ is the induced metric on $\Sigma_r$, the hypersurface with distance $r$ to $\Sigma_0$. It is not hard to see $\Sigma_r$ is convex and diffeomorphic to $\mathbf S^{n-1}$, so $g_r \in \mathcal M^k_{psc}(\mathbf S^{n-1}) $. Then 
\begin{equation*}
\gamma(t)=\left\{\begin{array}{cc}
(1-\log(1-t))^{-2}g_{-\log(1-t)}\quad\ \, t\in[0,1)\\
\gamma_{std}\qquad\qquad\qquad\qquad\qquad t=1
\end{array}\right.
\end{equation*}
is a continuous path in $\mathcal M^k_{psc}(\mathbf S^{n-1})$ joining $g_0$ and $\gamma_{std}$. Clearly, $\bar g$ is the standard Euclidean metric if we choose such $\gamma(t)$ in Lemma \ref{Lem: base space for Shi-Tam flow}. \end{exa}

For any smooth metric $\gamma$ on $\mathbf S^{n-1}$, we define
\begin{equation*}
\lambda_{min}(\gamma)=\sup\{\lambda>0\,|\,\gamma\geq \lambda\gamma_{std}\}
\end{equation*}
to measure the non-degeneracy of $\gamma$ with respect to the standard metric $\gamma_{std}$.
Consider the following class of metrics
\begin{equation*}
\mathcal M^n_{c_1,\,c_2,\,d,\,V}:=\left\{\gamma\in\mathcal M^\infty\left(\mathbf{S}^{n-1}\right)\bigg|
\begin{array}{c}
|\Rm_\gamma|\leq c_1,\,|\nabla_\gamma\Rm_\gamma|\leq c_2\\
\diam(\gamma)\leq d,\,\vol(\gamma)\geq V
\end{array}\right\}.
\end{equation*}
The following lemma states that for any metric in $\mathcal M^n_{c_1,\,c_2,\,d,\,V}$, one can always find a balanced parametrization, that is
\begin{lm}\label{Lem: good parametrization for M c_1 c_2 d}
There exists a universal constant $\Lambda=\Lambda(n,c_1,c_2,d,V)$ such that for any metric $\gamma\in\mathcal M^n_{c_1,\,c_2,\,d,\,V}$, we can find a diffeomorphism
$
\phi:\mathbf S^{n-1} \to\mathbf S^{n-1}
$
satisfying
\begin{equation*}
\|\phi^*\gamma\|_{C^2(\mathbf S^{n-1},\gamma_{std})}+\lambda_{min}^{-1}(\phi^*\gamma)\leq \Lambda.
\end{equation*}
\end{lm}

\begin{proof}
We argue by contradiction. Suppose the consequence is not true, then for any integer $k$, we can find a metric $\gamma_k\in\mathcal M^n_{c_1,\,c_2,\,d,\,V}$ such that for any diffeomorphism $\phi$, there holds
\begin{equation}\label{Eq: C2 norm unbounded}
\|\phi^*\gamma_k\|_{C^2(\mathbf S^{n-1},\gamma_{std})}+\lambda_{min}^{-1}(\phi^*\gamma_k)>k.
\end{equation}
However, it follows from the Cheeger-Gromov compactness theory that the space $\mathcal M^n_{c_1,\,c_2,\,d,\,V}$ is $C^{2,\alpha}$-precompact for any $0<\alpha<1$. Therefore, after passing to a subsequence (still denoted by $\gamma_k$), there exist diffeomorphisms
$
\phi_k:\mathbf S^{n-1}\to\mathbf S^{n-1}
$
such that $\tilde\gamma_k:=\phi_k^*\gamma_k$ converges to a limit metric $\tilde\gamma_\infty$ in the $C^{2,\alpha}$-sense (as metric functions in local coordinate charts). It is clear that the quantities
\begin{equation*}
\|\tilde\gamma_k\|_{C^2(\mathbf S^{n-1},\gamma_{std})}+\lambda_{min}^{-1}(\tilde\gamma_k)
\end{equation*}
converge to that of $\tilde\gamma_\infty$ under the $C^{2,\alpha}$-convergence, which contradicts (\ref{Eq: C2 norm unbounded}).
\end{proof}

\begin{lm}\label{Lem: path start from M c_1 c_2 d V_0}
For any metric $\gamma\in\mathcal M^n_{c_1,\,c_2,\,d,\,V}$, we can find a diffeomorphism $\phi:\mathbf S^{n-1} \to \mathbf S^{n-1}$ and a smooth path $\gamma_t:[0,1]\to \mathcal M^\infty(\mathbf S^{n-1})$ with $\gamma_0=\phi^*\gamma$ and $\gamma_t\equiv \gamma_{std}$ for $t\in[\frac56, 1]$ such that $|R_{\gamma_t}|$, $\|\gamma_t'\|_{\gamma_t}$ and $\|\gamma_t''\|_{\gamma_t}$ are bounded by universal constants depending only on $n$, $c_1$, $c_2$, $d$ and $V$.
\end{lm}
\begin{proof}
From Lemma \ref{Lem: good parametrization for M c_1 c_2 d}, there exists a diffeomorphism $\phi:\mathbf S^{n-1} \to \mathbf S^{n-1}$ such that
\begin{equation*}
\|\phi^*\gamma\|_{C^2(\mathbf S^{n-1},\gamma_{std})}+\lambda_{min}^{-1}(\phi^*\gamma)\leq \Lambda,
\end{equation*}
where $\Lambda$ is a universal constant depending only on $n$, $c_1$, $c_2$, $d$ and $V$. 

First we take a continuous path in $\mathcal M^\infty(\mathbf S^{n-1})$ by 
\begin{equation*}
\gamma^1_t=\left\{\begin{array}{cc}
(1-\frac{3}{2}t)\phi^*\gamma+\frac{3}{2}t\gamma_{std}\quad\  0\leq t\leq 2/3,\\
\gamma_{std}\qquad\qquad\qquad\qquad\ \, 2/3 < t\leq 1.
\end{array}\right.
\end{equation*}
Next we obtain a smooth path from $\gamma^1_t$ through a mollification procedure. Let $\varphi(t)$ be a smooth function with support in $(-1,1)$ that satisfies $0\leq\varphi\leq 1$, $\varphi(t)=\varphi(-t)$, and
\begin{equation*}
\int_{-\infty}^{+\infty}\varphi(t)\,\mathrm dt=1.
\end{equation*}
Let $\sigma$ be a fixed constant such that $0<\sigma\leq 1/6$ and $\varphi_{\sigma}(t)=\sigma^{-1}\varphi(\sigma^{-1}t)$. For $\frac{1}{2}\leq t\leq \frac{5}{6}$, we define
\begin{equation*}
\gamma^2_t=\varphi_\sigma\ast \gamma_t^1=\int_{-\sigma}^{\sigma}\varphi_{\sigma}(s)\gamma^1_{t-s}\,\mathrm ds.
\end{equation*}
Then it is not hard to see that the path
\begin{equation*}
\gamma_t=\left\{\begin{array}{cc}
\gamma_t^2&\frac{1}{2}\leq t\leq \frac{5}{6},\\
\gamma_t^1&\text{\rm elsewhere}
\end{array}\right.
\end{equation*}
is smooth and satisfies $\gamma_0=\phi^*\gamma$, $\gamma_t\equiv \gamma_{std}$ for $t\in[\frac56, 1]$. 

Thus, once we have proved $|R_{\gamma_t}|$, $\|\gamma_t'\|_{\gamma_t}$ and $\|\gamma_t''\|_{\gamma_t}$ are bounded by universal constants depending only on $\Lambda$, we reach our goal. Since $\gamma_t$ are convex combinations of $\phi^*\gamma$ and $\gamma_{std}$, we have the estimates $\|\gamma_t\|_{C^2(\mathbf S^{n-1},\gamma_{std})}\leq \Lambda$ and
\begin{equation*}
\min\{\Lambda^{-1},1\}\gamma_{std}\leq \gamma_t\leq \max\{\Lambda,1\}\gamma_{std}.
\end{equation*}
Therefore, $|R_{\gamma_t}|$ is bounded by a universal constant depending only on $\Lambda$. Note that the derivatives
\begin{equation*}
\gamma'_t=\left\{\begin{array}{cc}
\frac{3}{2}(\gamma_{std}-\phi^*\gamma)&0\leq t<1/2,\\
\varphi'_\sigma\ast \gamma_{t}^1&1/2\leq t\leq 5/6,\\
0& 5/6< t\leq 1
\end{array}\right.
\end{equation*}
and
\begin{equation*}
\gamma''_t=\left\{\begin{array}{cc}
\varphi''_\sigma\ast \gamma_{t}^1\  &1/2\leq t\leq 5/6,\\
0& \text{\rm elsewhere}
\end{array}\right.
\end{equation*}
are linear combinations of $\phi^*\gamma$ and $\gamma_{std}$, we conclude that the quantities $\|\gamma_t'\|_{\gamma_t}$ and $\|\gamma_t''\|_{\gamma_t}$ are also bounded by universal constants depending only on $\Lambda$.
\end{proof}

\begin{lm}\label{Lem: path M c_1 d}
Given any metric $\gamma\in\mathcal M^n_{c,\,d}$, we can find a diffeomorphism $\phi$ and a piecewise smooth path $\gamma_t:[0,1]\to \mathcal M^\infty(\mathbf S^{n-1})$ with $\gamma_0=\phi^*\gamma$ and $\gamma_t\equiv \gamma_{std}$ for $t\in[\frac56, 1]$ such that $R_{\gamma_t}$, $\|\gamma_t'\|_{\gamma_t}$ and $\|\gamma_t''\|_{\gamma_t}$ are bounded by universal constants depending only on $n$, $c$ and $d$ away from the unique broken point $t=1/3$. Furthermore, the path $\gamma_t$ is smooth on both sides of $t=1/3$ and satisfies $\gamma'_{(1/3)^+}=\gamma'_{(1/3)^-}=0$. 
\end{lm}
\begin{proof}
Let $\{\gamma^1_t\}_{0\leq t<T_s}$ be the Ricci flow with initial metric $\gamma$, where $T_s$ is the first singular time. It is standard that $\gamma^1_t$ is a smooth path in $\mathcal M^\infty(\mathbf S^{n-1})$. By Theorem 3.2.11 and Theorem 5.3.1 in \cite{PTopping}, we can find a universal positive constant $T=T(n,c)<T_s$ such that
\begin{equation}\label{Eq: curvature tensor estimate}
\|\Rm_t\|_{\gamma^1_t}\leq 2c\quad\  \text{for all}\ \, t\in[0, T].
\end{equation}
Having above estimates, by Theorem 3.3.1 in \cite{PTopping}, for any positive integer $k$ we get
\begin{equation}\label{Eq: curvature tensor derivative estimate}
\left\|\nabla^k_t\Rm_t\right\|_{\gamma^1_t}\leq C(n,k,c)t^{-\frac{k}{2}} \quad\  \text{for all}\ \, t\in\Big(0,\frac{1}{2c}\Big].
\end{equation}
Here and in the sequel, let $C(\cdot)$ denote universal constants depending only on quantities in the bracket. We may assume $T\leq\frac{1}{2c}$. For the path $\{\gamma^1_t\}_{t\in[0,T]}$, from the estimate (\ref{Eq: curvature tensor estimate}), we know that $|R_{\gamma^1_t}|$ and $\|(\gamma^1_t)'\|_{\gamma^1_t}$ are bounded by universal constants depending only on $n$ and $c$. Furthermore, it follows from the estimates (\ref{Eq: curvature tensor estimate}), (\ref{Eq: curvature tensor derivative estimate}) and the evolution equation
\begin{equation*}
\frac{\partial}{\partial t}\Rm_t=\Delta_t \Rm_t+\Rm_t\ast\Rm_t
\end{equation*}
that $\|(\gamma^1_t)''\|_{\gamma^1_t}\leq C(n,c)t^{-1}$ for all $t\in(0,T]$.
Define
\begin{equation*}
\hat\gamma^1_t:=\gamma^1_{c(t)}:[0,\sqrt[3]{6T}]\to \mathcal M^\infty(\mathbf S^{n-1}),
\end{equation*}
where $c(t)=\frac{\sqrt[3]{6T}}{2}t^2-\frac{1}{3}t^3$. Then the quantities $|R_{\hat\gamma^1_t}|$, $\|(\hat\gamma^1_t)'\|_{\hat\gamma^1_t}$ and $\|(\hat\gamma^1_t)''\|_{\hat\gamma^1_t}$ are bounded by universal constants depending only on $n$ and $c$, and $(\hat\gamma^1_{t})'\big|_{t=\sqrt[3]{6T}}=0$.

Next we construct another smooth path from the metric $\gamma^1_T$. From above discussion, one has
$\|\nabla_T\Rm_T\|_{\gamma^1_T}\leq C(n,c)$. It follows from the evolution equation
\begin{equation*}
{\partial_t}{\gamma^1_t}=-2\ric_{\gamma^1_t}
\end{equation*}
and the estimate (\ref{Eq: curvature tensor estimate}) that
\begin{equation*}
\diam(\mathbf S^{n-1},\gamma^1_T)\leq\diam(\mathbf S^{n-1},\gamma)e^{C(n,c)T}\leq C(n,c,d),
\end{equation*}
and
\begin{equation*}
\vol(\mathbf S^{n-1},\gamma^1_T)\geq \vol(\mathbf S^{n-1},\gamma_{std})e^{-C(n,c)T}\geq C(n,c)>0.
\end{equation*}
By Lemma \ref{Lem: path start from M c_1 c_2 d V_0}, we can find a diffeomorphism $\phi:\mathbf S^{n-1}\to \mathbf S^{n-1}$ and a smooth path $\gamma^2_t:[0,1]\to \mathcal M^\infty(\mathbf S^{n-1})$ with $\gamma^2_0=\phi^*(\gamma^1_T)$ and $\gamma^2_t\equiv \gamma_{std}$ for $t\in [\frac56, 1]$ such that quantities $|R_{\gamma^2_t}|$, $\|(\gamma^2_t)'\|_{\gamma^2_t}$ and $\|(\gamma^2_t)''\|_{\gamma^2_t}$ are bounded by universal constants depending only on $n$, $c$ and $d$. Actually we can further require $\gamma^2_t\equiv \phi^*(\gamma^1_T)$ 
around $t=0$ and hence $(\gamma^2_t)'\big |_{t=0}=0$.

Define a new path $\{\gamma_t\}_{t\in[0,1]}$ by
\begin{equation*}
\gamma_t=\left\{\begin{array}{cc}
\phi^*(\hat\gamma^1_{3\sqrt[3]{6T}t})&0\leq t\leq 1/3,\\
\gamma^2_{3t-1} &1/3<t\leq 2/3,\\
\gamma_{std} &2/3<t\leq 1.
\end{array}\right.
\end{equation*}
It is not hard to verify that the path $\{\gamma_t\}_{t\in[0,1]}$ satisfies all our requirements.
\end{proof}

The following lemma due to Brendle-Marques-Neves \cite{BMN} is very useful in gluing constructions.

\begin{lm}[Theorem 5 in \cite{BMN}]\label{BMNgluing}
Let $M$ be a compact manifold of dimension $n$ with boundary $\partial M$,
and let $g$ and $\tilde g$ be two smooth Riemannian metrics on $M$ such that
$g-\tilde g=0$ at each point on $\partial M$. Moreover, we assume that
$H_g-H_{\tilde g}> 0$ at each point on $\partial M$. Given any real number
$\varepsilon>0$ and any neighborhood $U$ of $\partial M$, there exists a
smooth metric $\hat g$ on M with the following properties:\\
$\bullet$ We have the pointwise inequality $R_{\tilde
g}(x)\geq\min\{R_g(x),R_{\tilde g}(x)\}-\varepsilon$ at each point $x\in
M$.\\
$\bullet$ $\hat g$ agrees with $g$ outside $U$.\\
$\bullet$ $\hat g$ agrees with $\tilde g$ in a neighborhood of $\partial M$.
\end{lm}

The following elementary lemma is used in the proof of Theorem \ref{Thm: Exponential Decay sigma Invariant}.
\begin{lm}\label{Lem: root of mu equation}
For any $\mu> 0$, there exists a unique root $c_\mu\in(0,1)$ of the equation
\begin{equation*}
x^{1-\frac{2}{n}}=\mu(1-x).
\end{equation*}
Furthermore, $c_\mu$ is a strictly monotone increasing continuous function of
$\mu$ with $\lim\limits_{\mu\to 0}c_\mu=0$.
\end{lm}

\begin{proof}
Define $f:(0,1)\rightarrow\mathbf{R}$ by
\begin{equation*}
f(x)=\frac{x^{1-\frac{2}{n}}}{1-x}.
\end{equation*}
Clearly, $f$ is smooth. It suffices to prove that
$f$ is strictly monotonically increasing with range $(0,+\infty)$. Through a
direct calculation, we get
\begin{equation*}
f'(x)=\frac{(1-\frac{2}{n})x^{-\frac{2}{n}}(1-x)+x^{1-\frac{2}{n}}}{(1-x)^2}>0,
\end{equation*}
and
\begin{equation*}
\lim_{x\to 0^+}f(x)=0,\quad \lim_{x\to 1^-}f(x)=+\infty.
\end{equation*}
\end{proof}

Next, we will construct the so called Schwarzschild neck, which is a PSC fill-in of a pair of Bartnik data $(\Sigma,\gamma,H_1)$ and $(\Sigma,\mu\gamma, H_2)$ with $H_1,\mu,H_2>0$. By virtue of this neck, we can extend the extent of the prescribed mean curvature that admits a PSC fill-in from $0$ to a certain positive constant. 
\begin{lm}\label{1}
Let $(\Sigma^{n-1},\gamma,H)$ be a triple of Bartnik data. Assume $R_{\gamma}>(n-1)(n-2)$ and $H$ is a constant in $(0,n-1)$. Let $h$ be a constant in $[0,H)$. Then the metric
 \[
 g=\psi^2(r)\left(\mathrm dr^2+r^2\gamma\right)
 \]
on $\Sigma\times[r_1,r_2]$ has the following properties:
\begin{enumerate}
\item[(1)] {$r_1\psi(r_1)<1$}, $r_2\psi(r_2)=1$, 
\item[(2)] $H_g|_{\Sigma\times\{r_1\}}\equiv h$, $H_g|_{\Sigma\times\{r_2\}}\equiv H$, 
\item[(3)] $R_g>0$, 
\end{enumerate}
where
$$ \psi(r)=\left(1+\frac{m}{2r^{n-2}}\right)^{\frac{2}{n-2}},\ \ \  m=\frac{1}{2}-\frac{H^{2}}{2(n-1)^{2}},$$
$${r_{1}}=r_1(h)\leq\left(\frac{m}{2}\right)^{\frac{1}{n-2}}\ \ \,\text{and}\ \ \  \,{r_{2}}=\left(\frac{n-1+H}{2(n-1)}\right)^{\frac{2}{n-2}}.$$
Furthermore, $(\Sigma\times[r_1,r_2],g)$ has an extension with positive scalar curvature.
 \end{lm}
 
\begin{proof}
Through a direct calculation, we get
\begin{equation*}
H_r=\frac{n-1}{r}\psi^{-\frac{n}{2}}(r)\left(1-\frac{m}{2r^{n-2}}\right),
\end{equation*}
where $H_r$ is the mean curvature of $\Sigma\times\{r\}$ with respect to the  $\partial_r$-direction, and
\begin{equation*}
R_g=r^{-2}\psi^{-2}\left(R_\gamma-(n-1)(n-2)\right).
\end{equation*}
Obviously, $R_g>0$. With the values of $m$ and $r_2$ given above, it is not hard to verify $r_2\psi(r_2)=1$ and $H_{r_2}=H$. Since $0\leq h<H$, we can solve $H_{r_1}=-h$ to get a unique root $r_1\in (0,(\frac{m}{2})^{\frac{1}{n-2}}]$.  Note that in our convention, $H_g|_{\Sigma\times\{r_1\}}=-H_{r_1}$ and $H_g|_{\Sigma\times\{r_2\}}=H_{r_2}$. Thus we obtain the desired manifold $(\Sigma\times[r_1,r_2],g)$. Choosing some
$r_1'<r_1$ and $r_2'>r_2$, then $(\Sigma\times[r_1',r_2'],g)$ gives the
extension.
\end{proof}

By scaling, the following result holds immediately.

\begin{pro}\label{neck}
Let $(\Sigma^{n-1},\gamma,H)$ be a triple of Bartnik data. Suppose $H$ is a positive constant and $R_\gamma>\frac{n-2}{n-1}H^2$.
Let $h$ be a constant in $[0,H)$. Then for any constant $\varepsilon$ satisfying $0<{\varepsilon}<\min R_\gamma-\frac{n-2}{n-1}H^2$, there exist a positive constant
$\mu<1$ and a metric $g_{\varepsilon}$, such that
$(\Sigma^{n-1},\gamma,H)$ and $(\Sigma^{n-1},\mu\gamma,h)$ can be realized as the
boundary data of the manifold $(\Sigma\times[r_1,r_2],g_{{\varepsilon}})$ with
$$R_{g_{\varepsilon}}\geq 
\min R_\gamma-\frac{n-2}{n-1}H^2-\varepsilon.$$
Furthermore, $(\Sigma^{n-1}\times[r_1,r_2] ,g_{\varepsilon})$ has an extension
with scalar curvature satisfying above inequality.
\end{pro}

Now, we give the definition of Schwarzschild neck.

\begin{defn}\label{Schwarzschild-neck}
Let $(\Sigma^{n-1},\gamma,H)$ be a triple of Bartnik data. Assume $H$ is a positive constant and $R_\gamma>\frac{n-2}{n-1}H^2$.
Then we call $(\Sigma^{n-1}\times[r_1,r_2],g_\epsilon)$ constructed in Proposition \ref{neck} a Schwarzschild neck of data $(\Sigma^{n-1},\gamma,H,h,\varepsilon)$.
\end{defn}

\section{Proof of Main Theorems} 
In this section, we prove our main results. 
\subsection{Non-existence of fill-in with NNSC metrics}\ 
\vskip 0.1cm
In this subsection, we give proofs of the results on non-existence of fill-in with NNSC metrics stated in the introduction, we first have:

\begin{proof}[Proof of Theorem \ref{nonexistfillin1}] It suffices to show there is some $h_0 < \infty$ satisfying Theorem \ref{nonexistfillin1}, as once we verify this fact we may take infimum to get the smallest one, which depends only on $\gamma$. Let $l=[\frac{k-1}2]\geq 2$. Fixing a small positive constant $\epsilon$, we can construct a $C^{2l+1}$ metric $\bar g$ on $\mathbf S^{n-1}\times[1,+\infty)$ as in Lemma \ref{Lem: base space for Shi-Tam flow}.
With the same notations in Lemma \ref{Lem: base space for Shi-Tam flow}, we consider the quasi-spherical metric equation
\begin{equation}\label{Eq: quasi spherical equation 2}
\left\{
\begin{aligned}
\bar H_s\frac{\partial u}{\partial s}&=u^2\Delta_{\bar\gamma_s}u+\frac{1}{2}(u-u^3)R_{\bar\gamma_s}-\frac{1}{2}R_{\bar g}u\\
\quad u(1)&=u_1>0,
\end{aligned}
\right.
\end{equation}
 where $\bar H_s={\mathrm t\mathrm r}_{\bar\gamma_s}\bar A_s$ is the mean curvature of $\Sigma_s=\mathbf S^{n-1}\times\{s\}$ with respect to $\bar g$ and the $\partial_s$-direction,  $\bar A_s$ is  the second fundamental form of $\Sigma_s$ with respect to the same direction, and $u_1$ is a smooth positive function on $\mathbf{S}^{n-1}$ to be given. Recall that $\tilde \gamma_s$ is a reparametrization of $\gamma_t$ and $\gamma_t$ is a smooth path in $\mathcal M_{psc}^{2l+1}(\mathbf S^{n-1})$, we see  that $R_{\bar\gamma_s}>0$ for all $s\geq 1$. Combined with the bounds on $R_{\bar\gamma_s}$ and $R_{\bar g}$, it follows from the parabolic maximum principle that the solution $u$ is positive and has bounded $C^0$ a priori estimate on any finite time interval. Therefore, equation (\ref{Eq: quasi spherical equation 2}) has a unique positive solution on the entire $[1,\infty)$. From the parabolic $L^p$-estimate and Schauder estimate, we conclude that $u$ is actually in H\"older space $C^{2l+\alpha,\,l+\alpha/2}$ for any $0<\alpha<1$.

Set 
$
g=u^2\mathrm ds^2+s^2\tilde\gamma_s,
$
then $g$ is a $C^{l,\,\alpha/2}$ metric. Let $A_s$ and $H_s$ denote the second fundamental and the mean curvature of $\Sigma_s$ induced from metric $g$. It is not hard to see 
\begin{equation}\label{meancurrelation}
A_s=u^{-1}\bar A_s,\qquad H_s=u^{-1}\bar H_s.
\end{equation}
By the Riccati equation, Gauss equation and relation \eqref{meancurrelation}, we have
\begin{equation*}
\begin{split}
\frac{\mathrm d}{\mathrm ds}\int_{\Sigma_s}H_s\,\mathrm d\mu_{\bar\gamma_s}&=\frac{1}{2}\int_{\Sigma_s}\left(\bar H_s^2-\|\bar A_s\|^2\right)u^{-1}\,\mathrm d\mu_{\bar\gamma_s}+\frac{1}{2}\int_{\Sigma_s}R_{\bar\gamma_s}u\,\mathrm d\mu_{\bar\gamma_s}\\
&\geq\frac{1}{2}\int_{\Sigma_s}\left(\bar H_s^2-\|\bar A_s\|^2\right)u^{-1}\,\mathrm d\mu_{\bar\gamma_s},
\end{split}
\end{equation*}
where we drop the second integral with the fact $R_{\bar\gamma_s}\geq 0$  in the second line.
Using estimate (\ref{Eq: mean curvature estimates on base space}) and relation \eqref{meancurrelation}, we see
\begin{equation*}
\begin{split}
\int_{\Sigma_s}\left(\bar H_s^2-\|\bar A_s\|^2\right)u^{-1}\,\mathrm d\mu_{\bar\gamma_s}&\geq \frac{(n-2)(1-\epsilon)}{s}\int_{\Sigma_s}\bar H_su^{-1}\,\mathrm d\mu_{\bar\gamma_s}\\
&=\frac{(n-2)(1-\epsilon)}{s}\int_{\Sigma_s}H_s\,\mathrm d\mu_{\bar\gamma_s}.
\end{split}
\end{equation*}
For convenience, we set
\begin{equation*}
\alpha(n,\epsilon)=\frac{(n-2)(1-\epsilon)}{2}.
\end{equation*}
Then we arrive at
\begin{equation*}
\frac{\mathrm d}{\mathrm ds}\int_{\Sigma_s}H_s\,\mathrm d\mu_{\bar\gamma_s}\geq\frac{\alpha(n,\epsilon)}{s}\int_{\Sigma_s}H_s\,\mathrm d\mu_{\bar\gamma_s}.
\end{equation*}
Integrating above ordinary differential inequality, we finally obtain 
\begin{equation*}
\int_{\Sigma_s}H_s\,\mathrm d\mu_{\bar\gamma_s} \geq s^{\alpha(n,\epsilon)}\int_{\Sigma_1}H_1\,\mathrm d\mu_{\bar\gamma_1}=s^{\alpha(n,\epsilon)}\int_{\mathbf S^{n-1}}\bar H_1u^{-1}_1\,\mathrm d\mu_{\gamma}.
\end{equation*}
Since $(\mathbf S^{n-1}\times[s_0,+\infty),\bar g)$ is Euclidean, it follows from Theorem 2.1 and Lemma 4.2 in \cite{Shi-Tam1} that $(\mathbf S^{n-1}\times[1,+\infty),g)$ is a scalar-flat AF end with ADM mass
\begin{equation}\label{mass}
\begin{split}
m_{ADM}&\leq C(n)\int_{\Sigma_{s_0}}\left(\bar H_{s_0}- H_{s_0}\right)\mathrm d\mu_{s_0}\\
&\leq C(n)\left(n(n-1)\omega_n s_{0}^{n-2}-s_0^{\alpha(n,\epsilon)}\int_{\mathbf S^{n-1}}\bar H_1u^{-1}_1\,\mathrm d\mu_{\gamma}\right),
\end{split}
\end{equation}
where $\omega_n$ is the volume of the unit ball in $\mathbf R^n$.

We now claim that the constant
\begin{equation*}
h_0=n(n-1)\omega_ns_0^{n-\alpha(n,\epsilon)-2}
\end{equation*}
satisfies our requirement. We argue by contradiction. Let $(\mathbf S^{n-1},\gamma,H)$ be a triple of Bartnik data with
\begin{equation*}
H>0,\quad \int_{\mathbf S^{n-1}}H\,\mathrm d\mu_\gamma>h_0,
\end{equation*}
and $(\Omega,\tilde g)$ be a fill-in of it with nonnegative scalar curvature. Setting $u_1=\bar H_1/H$, from above discussion we can obtain a scalar-flat AF end with negative ADM mass and inner boundary $(\mathbf S^{n-1},\gamma,H)$. Gluing $(\Omega,\tilde g)$ to this AF end, we get a complete AF manifold with NNSC and corners (see Definition 1 in \cite{Miao1}) along a closed hypersurface. Moreover, the mean curvatures on the two sides of the hypersurface are equal. By Theorem 1 in \cite{Miao1}, the ADM mass is nonegative. Thus we get the desired contradiction. 
\end{proof}

Next, we give:

\begin{proof}[Proof of Theorem \ref{Thm: no fillin nonnegative scalar 2}]
First we show that there exists a universal positive constant $H_0=H_0(n,c,d)$ such that for any $\gamma\in\mathcal M^n_{c,\,d}$, we can find a diffeomorphism $\phi:\mathbf S^{n-1}\to \mathbf S^{n-1}$ and a scalar-flat AF end $(E,g)$ that admits corners and has negative ADM mass and inner boundary $(\mathbf S^{n-1},\phi^*\gamma, H_0)$.

Given any $\gamma\in\mathcal M^n_{c,\,d}$, we can take a diffeomorphism $\phi:\mathbf S^{n-1}\to \mathbf S^{n-1}$ and a piecewise smooth path $\gamma_t:[0,1]\to \mathcal M^\infty(\mathbf S^{n-1})$ as in Lemma \ref{Lem: path M c_1 d}. Fixing $\epsilon$ to be a small positive constant, through a similar argument as in the proof of Lemma \ref{Lem: base space for Shi-Tam flow}, we can find a piecewise smooth metric $\bar g$ on $\mathbf S^{n-1}\times [1,+\infty)$ admitting corners along $\Sigma_{s_1}=\mathbf S^{n-1}\times\{s_1\}$ for some $1<s_1<s_0$. The second fundamental forms on two sides of $\Sigma_{s_1}$ are equal, and the estimates in Lemma \ref{Lem: base space for Shi-Tam flow} are still valid on each smooth piece. With $H_0$ to be determined later, using the notations in Lemma \ref{Lem: base space for Shi-Tam flow}, we consider the equation
\begin{equation}\label{Eq: quasi spherical equation 3}
\left\{
\begin{aligned}
\bar H_s\frac{\partial u}{\partial s}&=u^2\Delta_{\bar\gamma_s}u+\frac{1}{2}(u-u^3)R_{\bar\gamma_s}-\frac{1}{2}R_{\bar g}u\\
u(1)&=\frac{\bar H_1}{H_0}>0.
\end{aligned}
\right.
\end{equation}

Above equation holds on smooth parts of $(\mathbf S^{n-1}\times[1,+\infty),\bar g)$ and continuously cross the corners $\Sigma_{s_1}$. Note that $\bar g$ is exactly the Euclidean metric outside a compact set, for the long time existence for the solution $u$, we only need to rule out the possibility that $u$ blows up in a specific finite time interval. For this purpose, using bounds for $R_{\bar \gamma_s}$ and $R_{\bar g}$ in Lemma \ref{Lem: base space for Shi-Tam flow}, we can construct appropriate barrier functions from the corresponding ordinary differential equation of (\ref{Eq: quasi spherical equation 3}). With a comparison argument, we can take $H_0$ large enough, depending only on $n$, $c$ and $d$, such that the solution $u$ of (\ref{Eq: quasi spherical equation 3}) exists for all time and satisfies $0<u<1$. We also emphasize that $u$ is smooth on both sides of the corners $\Sigma_{s_1}$. Let $E=\mathbf S^{n-1}\times[1,+\infty)$ and $g=u^2\mathrm ds^2+s^2\tilde \gamma_s$. It follows from \cite{Shi-Tam1} that $(E,g)$ is the desired scalar-flat AF end, which admits corners and has negative mass and inner boundary $(\mathbf S^{n-1},\phi^*\gamma, H_0)$. On the two sides of the corners $\Sigma_{s_1}$, the second fundamental forms with respect to $g$ are equal. 

We now claim that $(\mathbf S^{n-1},\gamma,H)$ does not admit a fill-in with NNSC for any $H>H_0$. Otherwise, let $(\Omega,\tilde g)$ be one of such fill-ins. By gluing $(\Omega,\tilde g)$ and $(E,g)$ with the identification
$$\phi^{-1}:(\mathbf S^{n-1},\gamma)\to(\mathbf S^{n-1},\phi^*\gamma),$$
we obtain a complete AF manifold with NNSC but negative ADM mass. This AF manifolds has corners along two closed disjoint hypersurfaces, but the mean curvatures from both sides of the two hypersurfaces are equal. Notice that Theorem 1 in \cite{Miao1} is in fact valid for finite disjoint corners. Hence we obtain a contradiction.
\end{proof}

\subsection{Properties of $\theta$-invariant}\ 
\vskip 0.1cm
In this subsection, we prove the main results of $\theta$-invariant stated in the introduction. Let us introduce the following conventions first. 

When we use a symbol, for instance $\Omega$, to denote a fill-in region, we mean $\Omega$ denote the region and its boundary,  namely $\overline\Omega=\Omega$. And we use $\mathring{\Omega}$ to denote the interior of $\Omega$. If $(\Omega,g,X)$ is a fill-in of $(\Sigma,\gamma,H)$, by definition, we have $X^*(g|_{\partial\Omega})=\gamma$ and $H_g=H\circ X$. But in the following, for convenience, we omit $X$, just write $g|_{\partial\Omega}=\gamma$ and $H_g=H$.
For $n\geq 3$, let $c_n$ denote the dimensional constant $\frac{4(n-1)}{n-2}$. Let $C$ denote positive uniform constants in different situations with different values. 
\begin{proof}[Proof of Theorem \ref{staticity}]
We argue by contradiction. Suppose $(\Omega^{n},g,X)$ is an extremal fill-in of $\theta(\Sigma^{n-1},\gamma,H)$,
but $(\Omega^{n},g)$ is not static. Denote $\theta(\Sigma^{n-1},\gamma,H)$ by
$S$. By assumption, $g|_{\partial\Omega}=\gamma$, $H_g=H$ on $\partial\Omega$ and $R_g\geq S\geq 0$ in $\Omega$. We show the proof in three steps.

\emph{Step 1: Perturbation.}

Since $(\Omega^{n},g)$ is not static, according to Theorem 1 in \cite{Cor}, we can get a perturbed metric $g_1$
from $g$ that satisfies $g_1=g$ in a neighborhood of $\partial\Omega$, 
$R_{g_1}\geq R_g$ in $\Omega$, and $R_{g_1}(p)>R_g(p)$ for some $p\in\mathring\Omega$. As
$R_{g_1}(p)>R_g(p)\geq 0$, we can find a neighborhood $U_p$ of $p$, where
$R_{g_1}\geq(1+\rho)R_g$ for some constant $\rho>0$.

\emph{Step 2: Conformal deformation.}

We make conformal deformations to get a new metric $g_3$ that satisfies
$R_{g_3}>S$ in $\Omega$, $g_3=g$ and $H_{g_3}>H$ on $\partial\Omega$. We
discuss the following two cases.

\textbf{Case 1:} $S>0$.

In this case, we have to do twice conformal deformations.

Let $\eta$ be a smooth function compactly supported in $U_p$ that satisfies
$0\leq\eta\leq\frac{\rho}{1+\rho}$ and $\eta(p)=\frac{\rho}{1+\rho}$. Define
$f=\eta R_{g_1}$. Consider the following equation
\begin{equation}\label{conformal equation1}
\left\{
\begin{aligned}
\Delta_{g_1}u_1-c^{-1}_nfu_1&=0\quad\text{in}\ \Omega,\\
u_1&=1\quad \text{on}\ \partial\Omega.
\end{aligned}
\right.
\end{equation}
Since $f\geq 0$, above equation has a smooth solution $u_1$. By the maximum
principle, $0<u_1<1$ in $\mathring\Omega$ and $\frac{\partial
u_1}{\partial\nu}|_{\partial\Omega}>0$, where $\nu$ is the outward unit normal with
respect to $g$.

Let $g_2=u_1^\frac{4}{n-2}g_1$. Then
\begin{align*}
R_{g_2}&=u_1^{-\frac{n+2}{n-2}}\left(R_{g_1}u_1-c_n\Delta_{g_1}u_1\right)\\
&=u_1^{-\frac{4}{n-2}}(1-\eta)R_{g_1}.
\end{align*}
If $x\in U_p$, then
$$R_{g_2}(x)\geq(1+\rho)^{-1}u_1^{-\frac{4}{n-2}}(x)R_{g_1}(x)>R_g(x).$$
If $x\in\mathring\Omega\setminus U_p$, then $\eta(x)=0$ and
$$R_{g_2}(x)>R_{g_1}(x)\geq R_g(x).$$
So $R_{g_2}>S$ everywhere in $\mathring\Omega$. And
\begin{align*}
H_{g_2}=H_{g_1}+\frac{c_n}{2}\frac{\partial
u_1}{\partial\nu}=H+\frac{c_n}{2}\frac{\partial u_1}{\partial\nu}.
\end{align*}
So $H_{g_2}>H$ everywhere on $\partial\Omega$.

Thus we obtain a metric $g_2$ that satisfies $R_{g_2}>S$ in $\mathring\Omega$, $g_2=g$ and $H_{g_2}>H$ on
$\partial\Omega$.
In the following, we have to modify $g_2$ near  $\partial \Omega$. The key point is to find a positive smooth function $u_2$ on $\Omega$ that satisfies $u_2=1$ on $\partial \Omega$, $u_2\leq 1$ and $\Delta_{g_2}u_2<0$ near $\partial \Omega$. To that end, let $d(x)$ denote the distance function  from $x$ to $\partial\Omega$  and $\Omega_{\delta}$ denote the $\delta$-collar neighborhood of $\partial \Omega$ in
$\Omega$ with respect to $g_2$. Since $\partial \Omega$ is smooth, for sufficienltly small
$\delta$, $d$ is smooth in $\Omega_\delta$. We may assume $|\Delta_{g_2} d|\leq
C_1$ on $\overline\Omega_{\delta}$, where $C_1$ is a positive constant 
depending only on $\Omega_\delta$ and $g_2$. Let $w=(1-\beta d)^\alpha-1$
with constants $\alpha$ and $\beta $ to be determined later. Direct
calculation shows
$$\nabla_{g_2}w=-\alpha\beta\left(1-\beta d\right)^{\alpha-1}\nabla_{g_2}d,$$
and
\begin{equation*}
\begin{split}
\Delta_{g_2}w&=\alpha(\alpha -1)\beta ^2(1-\beta d)^{\alpha-2}-\alpha\beta(1-\beta
d)^{\alpha-1}\Delta_{g_2}d\\
&\leq \alpha(\alpha -1)\beta ^2(1-\beta d)^{\alpha-2}+\alpha\beta (1-\beta
d)^{\alpha-1}C_1.
\end{split}
\end{equation*}
Taking $\beta=2C_1$, we get
\begin{equation*}
\Delta_{g_2}w\leq 2\alpha C_1^2(1-\beta d)^{\alpha-2}\left(2\alpha-1-\beta
d\right).
\end{equation*}
Choosing $\alpha=1/4$ and sufficiently small $\delta_1$, we can find a
positive constant $\epsilon$ such that $\Delta_{g_2}w\leq -\epsilon<0$ in
$\Omega_{\delta_1}$. It is also easy to see that 
\begin{equation*}
\frac{\partial w}{\partial \nu}\Big|_{\partial\Omega}=\alpha\beta>0.
\end{equation*}
Now, we extend $w$ to the whole $\Omega$ to obtain a smooth function $v$ that satisfies $v<0$ in $\Omega$.
We may assume $\|v\|_{C^2(\Omega,g_2)}\leq C_2$ for some constant $C_2$.

Define $u_2=1+sv$, where $s$ is a small positive constant to be determined. Make
the following conformal deformation
\begin{equation*}
g_3=u_2^{\frac{4}{n-2}}g_2.
\end{equation*}
The scalar curvature after the conformal deformation is
\begin{equation*}
R_{g_3}=u_2^{-\frac{n+2}{n-2}}\left(R_{g_2}u_2-c_ns\Delta_{g_2}v\right).
\end{equation*}

Since $R_{g_2}>S$ in $\mathring\Omega$, there exists a positive constant $\epsilon'$
such that $R_{g_2}\geq S+\epsilon'$ in $\Omega\setminus\Omega_{\delta_1}$.
Therefore, in $\Omega\setminus\Omega_{\delta_1}$, 
\begin{align*}
R_{g_3}&\geq
R_{g_2}-c_nsu_2^{-\frac{n+2}{n-2}}\left|\Delta_{g_2}
v\right|\geq S+\epsilon'-O(s).
\end{align*}
Choosing sufficiently small $s$, we have $R_{g_3}>S$ in
$\Omega\setminus\Omega_{\delta_1}$.
Note that $u_2\leq 1$ and $\Delta_{g_2} v\leq -\epsilon$ in $\Omega_{\delta_1}$. 
Therefore we get 
$$
R_{g_3}\geq R_{g_2}+c_ns\epsilon\geq S+c_ns\epsilon
$$
in $\Omega\setminus\Omega_{\delta_1}$. Consequently, $R_{g_3}>S$ in $\Omega$. On $\partial\Omega$, we have
\begin{equation*}
H_{g_3}=H_{g_2}+\frac{c_ns}{2}\frac{\partial v}{\partial\nu}>H.
\end{equation*}

\textbf{Case 2:} $S=0$.

In this case, we only need to do conformal deformation once. For
$\varepsilon\geq 0$, consider the following equation
\begin{equation*}
\left\{
\begin{aligned}
\Delta_{g_1}u_\varepsilon-c^{-1}_nfu_\varepsilon&=-\varepsilon\quad\,\,\text{in}\ 
\,\Omega,\\
u_\varepsilon&=1\qquad\, \text{on}\ \, \partial\Omega.
\end{aligned}
\right.
\end{equation*}
For sufficiently small $\varepsilon$, above equation has a positive
smooth (with respect to both variables and the parameter $\varepsilon$) solution
$u_{\varepsilon}$. Let $g_3=u_{\varepsilon}^\frac{4}{n-2}g_1$. Then
\begin{align*}
R_{g_3}&=u_{\varepsilon}^{-\frac{n+2}{n-2}}\left(R_{g_1}u_{\varepsilon}-c_n\Delta_{g_1}u_{\varepsilon}\right)\\
&=u_{\varepsilon}^{-\frac{n+2}{n-2}}\left((1-\eta)R_{g_1}u_{\varepsilon}+c_n\varepsilon
\right)\\
&\geq c_n\varepsilon u_{\varepsilon}^{-\frac{n+2}{n-2}}.
\end{align*}
So for $\varepsilon>0$,  $R_{g_3}>0$ in $\Omega$. On $\partial\Omega$, we have
\begin{align*}
H_{g_3}=H+\frac{c_n}{2}\frac{\partial u_\varepsilon}{\partial\nu},
\end{align*}
where $\nu$ is the outward unit normal with respect to $g$.
When $\varepsilon=0$, $u_0$ satisfies \eqref{conformal equation1}. Since
$\frac{\partial u_0}{\partial\nu}|_{\partial\Omega}>0$ and $u_{\varepsilon}$ depends
smoothly on $\varepsilon$, for sufficiently small $\varepsilon$, $H_{g_3}>H$.

\emph{Step 3: Gluing}.

Now, we will use a similar argument as in the proof of Lemma 20 in \cite{J} to complete the proof. Roughly speaking, we will construct a metric $g_4$ on a small collar neighborhood of $\partial\Omega$ that satisfies $R_{g_4}>S$ in this collar neighborhood, and $g_4=g$, $H_{g_4}=H$ on $\partial\Omega$. Then we 
glue $g_4$ and $g_3$ to get a new metric $g_5$ that satisfies $R_{g_5}>S$ in
$\Omega$, $g_5|_{\partial\Omega}=\gamma$, and $H_{g_5}=H$. Thus we get a
contradiction.

For some small $t_0>0$, $\Sigma\times[-t_0,0]$ is diffeomorphic to a $t_0$-collar
neighborhood of $\partial\Omega$ in $\Omega$ with respect to $g_3$. Let $\Sigma_t$ denote $\Sigma\times\{t\}$ and identify $\Sigma\times\{0\}$ with $\partial\Omega$.  In this $t_0$-collar
neighborhood, we write $g_3$ as
$g_3(t)=\mathrm dt^2+\hat g_3(t)$, where $\hat g_3(t)$ is the metric on
$\Sigma_t$ induced from $g_3$. Define $\omega:\Sigma\rightarrow\mathbb{R}$ by
$$\omega(y)=\frac{H_{g_3}(y)-H(y)}{n-1}.$$ By definition, $\omega>0$. Let $\kappa$ be a smooth function on $[-t_0,0]$
that satisfies $\kappa(0)=0$, $\kappa'(0)=-1$. In $\Sigma\times [-t_0,0]$, define
$$g_4(y,t)=\mathrm dt^2+\left(1+\omega(y)\kappa(t)\right)^2\hat g_3(y,t).$$
Then extend $g_4$ to the whole $\Omega$ (in an arbitrary manner). Obviously, $g_4 |_{\partial\Omega}=\gamma$ and $\frac{\partial}{\partial t}|_{\partial\Omega}=\nu$. Use $\hat g_4(t)$ to denote the metric on $\Sigma_t$ induced from $g_4$. 

Let $A_i(t)$ and $H_i(t)$ denote the second fundamental form and the mean curvature of $\Sigma_t$ with respect to $g_i$ and the $\partial_t$-direction, for $i=3,4$.
We have
$$A_4=\left(1+\omega\kappa\right)^2A_3+\omega\kappa'\left(1+\omega\kappa\right)\hat
g_3,$$ and
\begin{align*}
H_4=H_3+\frac{(n-1)\omega\kappa'}{1+\omega\kappa}.
\end{align*}
So $H_{g_4}=H_4(0)=H$.

Let $\hat R_i(t)$ denote the scalar curvature of $\Sigma_t$ with respect $\hat g_i$ ($i=3,4$). We have 
\begin{align*}
\hat R_4=\left(1+\omega\kappa\right)^{-2}\left(\hat R_3-\frac{2(n-1)\kappa}{1+\omega\kappa}\Delta_{\hat g_4}\omega-\frac{(n-1)(n-4)\kappa^2}{(1+\omega\kappa)^2}|\nabla_{\hat g_4}\omega|^2\right).\end{align*}
By the Riccati equation and Gauss equation,
$$R_{g_i}=-2\frac{\partial H_i}{\partial t}+\hat R_i-|H_i|^2-|A_i|^2\quad
(i=3,4).$$
So
\begin{align*}
R_{g_4}=&R_{g_3}-2\frac{\partial(H_4-H_3)}{\partial t}+(\hat R_4-\hat
R_3)-(|H_4|^2-|H_3|^2)-(|A_4|^2-|A_3|^2)\\
=&R_{g_3}-\frac{2(n-1)\omega\kappa''}{1+\omega\kappa}-\frac{(n-1)(n-2)\omega^2\kappa'^2}{(1+\omega\kappa)^2}-\frac{\omega\kappa(2+\omega\kappa)}{(1+\omega\kappa)^{2}}\hat
R_3\\
&-\frac{2(n-1)\kappa}{(1+\omega\kappa)^3}\Delta_{\hat g_4}\omega-\frac{(n-1)(n-4)\kappa^2}{(1+\omega\kappa)^4}|\nabla_{\hat g_4}\omega|^2-\frac{2n\omega\kappa'}{1+\omega\kappa}H_3.
\end{align*}

Note that $H_3$, $\hat R_3$, $\Delta_{\hat g_4}\omega$ and $|\nabla_{\hat g_4}\omega|^2$ are bounded in $[-t_0,0]$. If $\kappa''(t)\ll -1$ in a small interval
around $t=0$, then $R_{g_4}>S$ in this small interval. We assume the interval is $[-t_1,0]$, for some $t_1< t_0$.

According to {\it Step 2} and above paragraph, there exists a
$\epsilon_1>0$ such that $R_{g_3}\geq S+\epsilon_1$ in $\Omega$ and
$R_{g_4}\geq S+\epsilon_1$ in $\Sigma\times [-t_1,0]$. To
glue $g_4$ to $g_3$, we apply Lemma \ref{BMNgluing} to the setting: $M=\Omega$, $g=g_3$, $\tilde g=g_4$, $\varepsilon=\epsilon_1/2$
and $U=\Sigma\times [-t_1,0]$. Then we get a new metric
$g_5=\hat g$ on $\Omega$. By the third property in Lemma \ref{BMNgluing}, $g_5$ agrees with $g_4$ in a neighborhood of $\partial\Omega$, so $g_5|_{\partial\Omega}=\gamma$ and $H_{g_5}=H$. By the second property, when
$x\in\Omega\setminus U$, $g_5=g_3$, so $R_{g_5}(x)=R_{g_3}(x)\geq S+\epsilon_1$. When $x\in
U$, according to the first property in Lemma \ref{BMNgluing}, $R_{g_5}(x)\geq S+\epsilon_1/2$. Hence,
$R_{g_5}(x)\geq S+\epsilon_1/2$ for all  $x\in\Omega$. Consequently, we get a
contradiction.
\end{proof}

Next, we prove Theorem \ref{Thm: weak monotonicity sigma invariant}, which is on the monotonicity of the $\theta$-invariant with respect to the prescribed mean curvature. 
\begin{proof}[Proof of Theorem \ref{Thm: weak monotonicity sigma invariant}]
We take the contradiction argument. If $\theta(\Sigma,\gamma,H_1)>\theta(\Sigma,\gamma,H_2)$,
then there exist a fill-in 
$(\Omega,g)$ of $(\Sigma, \gamma , H_1)$ and a 
 positive constant $\epsilon$ such that
$R_g\geq\theta(\Sigma,\gamma,H_2)+\epsilon$ in $\Omega$. 
Then we make a conformal transform of $g$ to increase
$H_g$ but not decrease $R_g$ much. As in {\it Step 2} of the proof of Theorem \ref{staticity}, let
$$u=1+sv\quad\  \text{and} \quad\   g'=u^{\frac{4}{n-2}}g,$$
where $v$ is the function defined in {\it Step 2} in the proof of Theorem \ref{staticity} and $s$ is a positive constant to
be determined later. Under this conformal transformation, the scalar
curvature of $g'$ is
\begin{equation*}
R_{g'}=u^{-\frac{n+2}{n-2}}\left(R_gu-c_ns\Delta_{g}v\right).
\end{equation*}
Therefore,
\begin{align*}
R_{g'}&\geq
R_g-\left|\left(u^{-\frac{4}{n-2}}-1\right)R_g-c_nsu^{-\frac{n+2}{n-2}}\Delta_g
v\right|\\&\geq \theta(\Sigma,\gamma,H_2)+\epsilon-O(s).
\end{align*}
Choosing $s$ sufficiently small, we have $R_{g'}\geq
\theta(\Sigma,\gamma,H_2)+\epsilon/2$ in $\Omega$. The mean curvature of $\partial\Omega$ under this conformal
deformation is
\begin{equation*}
H_{g'}=H_g+\frac{c_ns}{2}\frac{\partial v}{\partial\nu}>H_1,
\end{equation*}
Carrying out {\it Step 3} in the proof of Theorem \ref{staticity}, we can get a new metric $\hat
g$ such that $R_{\hat g}>\theta(\Sigma,\gamma,H_2)$ in $\Omega$, $\hat g|_{\partial\Omega}=\gamma$ and
$H_{\hat g}=H_2$. This contradicts the definition of
$\theta(\Sigma,\gamma,H_2)$. Hence, $\theta(\Sigma,\gamma,H_1)\leq
\theta(\Sigma,\gamma,H_2)$.
\end{proof}

Now, we begin to prove Theorem \ref{Thm: Exponential Decay sigma Invariant}. Once Theorem \ref{Thm: Exponential Decay sigma Invariant} is proved, Corollary \ref{limofsigma} follows immediately. We start with the following proposition, which states that the
$\theta$-invariant will decrease for a certain portion when the mean
curvature is lifted from $1$ to some constant $\lambda>1$.
\begin{pro}\label{Prop: strict inequality 2}
Let $(\Sigma^{n-1},\gamma)$ be a Riemannian manifold with $R_{\gamma}\geq 0$. Given a constant $\lambda>1$, if
$\theta(\Sigma^{n-1},\gamma,\lambda)<+\infty$, then there exists a constant
$\alpha>1$ such that
\begin{equation}\label{Eq: strict inequality}
\theta(\Sigma,\gamma,\lambda)\leq \alpha^{-2}\theta(\Sigma,\gamma,1).
\end{equation}
Furthermore, when $\theta(\Sigma^{n-1},\gamma,\lambda)>0$, we can choose
\begin{equation*}
\alpha=(1-c_\mu)^{\frac{1}{n-2}}\left(\lambda^2+\frac{n-1}{n}\theta(\Sigma,
\gamma,\lambda)\right)^{\frac{1}{n-2}},
\end{equation*}
where $c_\mu$ is the unique positive solution to $x^{1-\frac{2}{n}}=\mu(1-x)$ and
\begin{equation}\label{muandsigma}
\mu=\frac{n-1}{n}\left(\frac{n}{n-1}\lambda^2+\theta(\Sigma,\gamma,\lambda)
\right)^{\frac{2}{n}}
\theta(\Sigma,\gamma,\lambda)^{1-\frac{2}{n}}.
\end{equation}
\end{pro}
\begin{proof}
If $\theta(\Sigma,\gamma,\lambda)\leq 0$, with the choice $\alpha=2$, inequality (\ref{Eq: strict inequality}) is actually trivial. That is, if
$\theta(\Sigma,\gamma,1)\geq 0$, we have
\begin{equation*}
\theta(\Sigma,\gamma,\lambda)\leq 0\leq \alpha^{-2}\theta(\Sigma,\gamma,1).
\end{equation*}
Otherwise, from the monotonicity of $\theta$-invariant, we see
\begin{equation*}
\theta(\Sigma,\gamma,\lambda)\leq \alpha^{-2}\theta(\Sigma,\gamma,\lambda)\leq
\alpha^{-2}\theta(\Sigma,\gamma,1).
\end{equation*}

In the following, we deal with the case
$\theta(\Sigma,\gamma,\lambda)>0$.
By the definition of $\theta$-invariant, for any $\epsilon>0$, we can find a
fill-in $(\Omega_\epsilon,g_\epsilon)$ of $(\Sigma,\gamma,\lambda)$ with
$R_{g_\epsilon}\geq \theta(\Sigma,\gamma,\lambda)-\epsilon$. Then we construct
a neck $(\Omega_{neck},g_{neck})$ with
$R_{g_{neck}}\geq\theta(\Sigma,\gamma,\lambda)$ to glue to
$(\Omega_\epsilon,g_\epsilon)$. The boundary data of $(\Omega_{neck},g_{neck})$ is
$(\Sigma,\gamma,\lambda-\epsilon)$ and
$(\Sigma,\alpha_\epsilon^2\gamma,\alpha_\epsilon^{-1})$, where $\alpha_\epsilon>1$ is a constant to be
determined later. In detail, consider the manifold $(\Sigma\times
(0,\frac{\pi}{n\sigma}],g_\sigma)$ with
\begin{equation*}
g_\sigma=\mathrm
dt^2+c^2\sigma^{-2}\sin^{\frac{4}{n}}\left(\frac{n\sigma}{2}t\right)\gamma,
\end{equation*}
where the constants $c$ and $\sigma$ will be given later.
For simplicity, let 
\begin{equation*}
\alpha(t)=c\sigma^{-1} \sin^{\frac{2}{n}}\left(\frac{n\sigma}{2}t\right).
\end{equation*}
It follows from a straightforward calculation that
\begin{equation*}
R_{g_\sigma}=n(n-1)\sigma^2+\alpha^{-2}R_\gamma\geq n(n-1)\sigma^2,
\end{equation*}
where we have used the fact $R_\gamma\geq 0$ in the second step. To guarantee
$R_{g_\sigma}\geq \theta(\Sigma,\gamma,\lambda)$, we choose $\sigma$ that satisfies
$n(n-1)\sigma^2=\theta(\Sigma,\gamma,\lambda)$. Let $H_t$ denote the mean
curvature of $\Sigma\times\{t\}$ with respect to the $\partial_t$-direction and $\tilde
H_t$ denote the normalized mean curvature $\alpha(t)H_t$. Then we have 
\begin{equation*}
H_t=(n-1)\sigma\cot\left(\frac{n\sigma}{2}t\right),
\end{equation*}
and
\begin{equation*}
\tilde
H_t=c(n-1)\cot\left(\frac{n\sigma}{2}t\right)\sin^{\frac{2}{n}}
\left(\frac{n\sigma}{2}t\right).
\end{equation*}
Choosing $t_1$ to guarantee $H_{t_1}=\lambda-\epsilon$, fixing $c$ to ensure
$\alpha(t_1)=1$, and selecting $t_2$ such that $\tilde H_{t_2}=1$, we now
define $(\Omega_{neck},g_{neck})=(\Sigma\times[t_1,t_2],g_\sigma)$. It is not
difficult to see that its boundary data is
$(\Sigma,\gamma,\lambda-\epsilon)$ and
$(\Sigma,\alpha_\epsilon^2\gamma,\alpha_\epsilon^{-1})$ with
$\alpha_\epsilon=\alpha(t_2)$. 

Choose $t'_1\in (0,t_1)$. Then $\Sigma\times[t'_1,t_1]$ is diffeomorphic to a collar
neighborhood of $\partial\Omega_\epsilon$ in $\Omega_\epsilon$ with $\Sigma\times\{t_1\}$ diffeomorphic to $\partial\Omega_\epsilon$. Under this diffeomorphism, $g_\sigma$ is defined on a neighborhood of $\partial\Omega_\epsilon$. We extend $g_{\sigma}$ smoothly to the whole $\Omega_\epsilon$ (in an arbitrary manner). Note $g_\epsilon=g_\sigma$ and $\lambda=H_{g_\epsilon}>H_{g_\sigma}=\lambda-\epsilon$ on $\partial\Omega_\epsilon$. Applying Lemma \ref{BMNgluing} to the setting: $M=\Omega_\epsilon$, $g=g_\epsilon$, $\tilde g=g_{\sigma}$, $\varepsilon=\epsilon$ and $U=\Sigma\times [t'_1,t_1]$, we can get a new metric $\hat g$ on $\Omega_\epsilon$ that satisfies $R_{\hat g}\geq\theta(\Sigma,\gamma,\lambda)-2\epsilon$ in $\Omega_\epsilon$ and $\hat g=g_{\sigma}$ in a neighborhood of $\partial\Omega_\epsilon$. Define $\bar\Omega=\Omega_\epsilon\coprod\Omega_{neck}/\sim$, where $\sim$ is the diffeomorphism between $\partial\Omega_\epsilon$ and $\Sigma\times\{t_1\}$. Then Define $\bar g$ on $\bar\Omega$ by 
\begin{equation*}
\bar g=\left\{
\begin{aligned}
&\hat g \qquad\ \ x\in\Omega_\epsilon,\\
&g_{neck}\quad\  x\in\Omega_{neck}.
\end{aligned}
\right.
\end{equation*}
It is not hard to see that $\bar g$ is smooth and $(\bar\Omega,\bar g)$ gives a fill-in of $(\Sigma,\alpha_\epsilon^2\gamma,\alpha_\epsilon^{-1})$ with $R_{\bar g}\geq\theta(\Sigma,\gamma,\lambda)-2\epsilon$. By rescaling, we see
\begin{equation*}
\theta(\Sigma,\gamma,1)=\alpha_\epsilon^2\theta(\Sigma,\alpha_\epsilon^2
\gamma,\alpha_\epsilon^{-1})\geq
\alpha_\epsilon^{2}(\theta(\Sigma,\gamma,\lambda)-2\epsilon).
\end{equation*}
Letting $\epsilon\to 0$, we obtain (\ref{Eq: strict inequality}) with
$\alpha=\lim\limits_{\epsilon\to 0}\alpha_\epsilon$.

In the following, we calculate the explicit value of $\alpha$. First we list the equations in our construction as following:
\begin{equation}\label{Eq: 1}
n(n-1)\sigma^2=\theta(\Sigma,\gamma,\lambda),
\end{equation}
\begin{equation}\label{Eq: 2}
H_{t_1}=(n-1)\sigma\cot\left(\frac{n\sigma}{2}t_1\right)=\lambda-\epsilon,
\end{equation}
\begin{equation}\label{Eq: 3}
\alpha(t_1)=c\sigma^{-1}\sin^{\frac{2}{n}}\left(\frac{n\sigma}{2}t_1\right)=1,
\end{equation}
\begin{equation}\label{Eq: 4}
\tilde H_{t_2}=c(n-1)\cot\left(\frac{n\sigma}{2}t_2\right)\sin^{\frac{2}{n}}
\left(\frac{n\sigma}{2}t_2\right)=1,
\end{equation}
\begin{equation}\label{Eq: 5}
\alpha_\epsilon=\alpha(t_2)=c\sigma^{-1}\sin^{\frac{2}{n}}\left(\frac{n\sigma}{2}t_2\right).
\end{equation}
Dividing (\ref{Eq: 5}) by (\ref{Eq: 3}), we see
\begin{equation*}
\alpha_\epsilon=\left(\frac{\sin\left(\frac{n\sigma}{2}t_2\right)}
{\sin\left(\frac{n\sigma}{2}t_1\right)}\right)^{\frac{2}{n}}.
\end{equation*}
Multiplying (\ref{Eq: 2}) by (\ref{Eq: 3}), then dividing the obtained equation by (\ref{Eq: 4}),
we get
\begin{equation}\label{Eq: 8}
\lambda-\epsilon=\frac{\cos\left(\frac{n\sigma}{2}t_1\right)}
{\cos\left(\frac{n\sigma}{2}t_2\right)}
\left(\frac{\sin\left(\frac{n\sigma}{2}t_2\right)}
{\sin\left(\frac{n\sigma}{2}t_1\right)}\right)^{1-\frac{2}{n}}
=\frac{\cos\left(\frac{n\sigma}{2}t_1\right)}
{\cos\left(\frac{n\sigma}{2}t_2\right)}\alpha_\epsilon ^{\frac{n}{2}-1}.
\end{equation}
From (\ref{Eq: 1}) and (\ref{Eq: 2}), we have
\begin{equation}\label{Eq: 6.9}
\cos^2\left(\frac{n\sigma}{2}t_1\right)=\frac{n(\lambda-\epsilon)^2}
{n(\lambda-\epsilon)^2+(n-1)\theta(\Sigma,\gamma,\lambda)}
\end{equation}
and
\begin{equation}\label{Eq: 7}
\sin^2\left(\frac{n\sigma}{2}t_1\right)=\frac{(n-1)\theta(\Sigma,\gamma,\lambda)}
{n(\lambda-\epsilon)^2+(n-1)\theta(\Sigma,\gamma,\lambda)}.
\end{equation}
Combining (\ref{Eq: 3}), (\ref{Eq: 4}) and (\ref{Eq: 7}) together, we obtain
\begin{equation}\label{Eq: 6}
\sin^{2-\frac{4}{n}}\left(\frac{n\sigma}{2}t_2\right)=\mu_\epsilon\left(1-\sin^{2}
\left(\frac{n\sigma}{2}t_2\right)\right),
\end{equation}
where
\begin{equation*}
\mu_\epsilon=\frac{n-1}{n}\left(\theta(\Sigma,\gamma,\lambda)+
\frac{n}{n-1}(\lambda-\epsilon)^2\right)^{\frac{2}{n}}
\theta(\Sigma,\gamma,\lambda)^{1-\frac{2}{n}}.
\end{equation*}
By Lemma \ref{Lem: root of mu equation}, we can solve (\ref{Eq: 6}) to
obtain
\begin{equation*}
\sin^{2}\left(\frac{n\sigma}{2}t_2\right)=c_{\mu_\epsilon}\in(0,1).
\end{equation*}
Combing above equation with \eqref{Eq: 8} and  (\ref{Eq: 6.9}), we arrive at 
\begin{equation*}
\alpha_\epsilon=(1-c_{\mu_\epsilon})^{\frac{1}{n-2}}
\left((\lambda-\epsilon)^2+\frac{n-1}{n}\theta(\Sigma,
\gamma,\lambda)\right)^{\frac{1}{n-2}}.
\end{equation*}
Letting $\epsilon\to 0$, by the continuity of $c_{\mu_\epsilon}$ from
Lemma \ref{Lem: root of mu equation}, we obtain the desired result.
\end{proof}

We fix $\lambda=2$ to obtain the following corollary:
\begin{cor}\label{Cor: Uniform Decay}
There exist dimensional constants $\theta_0>0$ and $\alpha_0>1$ such that for any Riemannian manifold $(\Sigma^{n-1},\gamma)$ with $R_{\gamma}\geq 0$ and $\theta(\Sigma^{n-1},\gamma,2)\leq\theta_0$, the following holds 
\begin{equation*}
\theta(\Sigma,\gamma,2)\leq \alpha_0^{-2}\theta(\Sigma,\gamma,1).
\end{equation*}
\end{cor}
\begin{proof}
By the monotonicity of $\theta$-invariant, we only need to deal with the case $\theta(\Sigma,\gamma,2)>0$. Fixing $\lambda=2$ in Proposition \ref{Prop: strict inequality 2}, we have
\begin{equation*}
\theta(\Sigma,\gamma,2)\leq\alpha^{-2}\theta(\Sigma,\gamma,1),
\end{equation*}
where
\begin{equation*}
\alpha=(1-c_\mu)^{\frac{1}{n-2}}\left(4+\frac{n-1}{n}\theta(\Sigma,
\gamma,2)\right)^{\frac{1}{n-2}}> [4(1-c_\mu)]^{\frac{1}{n-2}},
\end{equation*}
with $c_\mu$ as in Proposition \ref{Prop: strict inequality 2}. Regard $c_\mu$ as a function of $\theta(\Sigma,\gamma,2)$ by the relation \eqref{muandsigma} and note that
$c_\mu$ converges to $0$ as $\theta(\Sigma,\gamma,2)$ tends to 0. So there exists a $\theta_0>0$ such that
 for any $(\Sigma,\gamma)$ with $\theta(\Sigma,\gamma,2)\leq\theta_0$, $4(1-c_\mu)\geq 2$. Then we can take
$\alpha_0=2^{\frac{1}{n-2}}$.
\end{proof}

Now we prove Theorem \ref{Thm: Exponential Decay
sigma Invariant} by iteration. 
\begin{proof}[Proof of Theorem \ref{Thm: Exponential Decay sigma
Invariant}]
We may assume that there exists a positive constant $H_1$ such that $\theta(\Sigma,\gamma,H_1)<+\infty$. Otherwise, case (1) holds. If there is a constant $H_2\geq H_1$ such that $\theta(\Sigma,\gamma,H_2)\leq 0$, by the monotonicity of $\theta$-invariant, case (2) holds trivially with $H_0=H_2$. 
Otherwise,  $\theta(\Sigma,\gamma,H)>0$ for any constant $H\geq H_1$. 

For those $H$, by rescaling, we have 
\begin{equation}\label{Eq: Exponential Decay sigma Invariant 1}
\theta(\Sigma,\gamma,2H)=H^2\theta(\Sigma,H^2\gamma,2).
\end{equation}
Substituting \eqref{Eq: Exponential Decay sigma Invariant 1} in the relation
$$0<\theta(\Sigma,\gamma,2H)\leq\theta(\Sigma,\gamma,H_1)<+\infty,$$ we see
\begin{equation*}
\theta(\Sigma,H^2\gamma,2)=O(H^{-2})\quad \text{\rm as}\ \ H\to+\infty.
\end{equation*}
Therefore, we can find a constant $H_0\geq \max\{1,H_1\}$ so that
$\theta(\Sigma,H^2\gamma,2)\leq \theta_0$ for any $H\geq H_0$, where
$\theta_0$ is the constant in Corollary \ref{Cor: Uniform Decay}. According to 
Corollary \ref{Cor: Uniform Decay}, there exists an absolute constant
$\alpha_0>1$ such that
\begin{equation}\label{intermediate}
\theta(\Sigma,H^2\gamma,2)\leq \alpha_0^{-2}\theta(\Sigma,H^2\gamma,1).
\end{equation}
Combining (\ref{Eq: Exponential Decay sigma Invariant 1}) with \eqref{intermediate} and rescaling,
we obtain
\begin{equation*}
\theta(\Sigma,\gamma,2H)\leq
\alpha_0^{-2}\theta(\Sigma,\gamma,H)\quad\text{for }\, H\geq H_0.
\end{equation*}
By iteration, it is clear that
\begin{equation*}
\theta(\Sigma,\gamma,2^kH_0)\leq \alpha_0^{-2k}\theta(\Sigma,\gamma,H_0)\quad
\text{for}\ \,k\in\mathbf N.
\end{equation*}
For any $H\geq H_0$, there is a $k\in\mathbf N_+$ such that $2^{k-1}H_0\leq
H<2^kH_0$. From this, we deduce $k\geq \log_2(H/H_1)$ and further
\begin{align*}
\theta(\Sigma,\gamma,H)&\leq \theta(\Sigma,\gamma,2^{k-1}H_0)\\&\leq
\alpha_0^{-2(k-1)}\theta(\Sigma,\gamma,H_0)\\ &\leq
\alpha_0^2\theta(\Sigma,\gamma,H_0)H_0^{2\log_2\alpha_0}H^{-2\log_2\alpha_0}.
\end{align*}
Taking $\beta=2\log_2\alpha_0$ and $C=\alpha_0^2\theta(\Sigma,\gamma,H_0)H_0^{2\log_2\alpha_0}$, 
we obtain the desired decay estimate. 
\end{proof}

Before we prove Theorem \ref{nonnegative}, we establish two propositions, which are about two interesting properties of the $\theta$-invariant. 
\begin{pro}\label{neqRbdy}
If $\theta(\Sigma,\gamma,0)\neq 0$, then $\theta(\Sigma,\gamma,0)\geq \min R_\gamma$.
\end{pro}

\begin{proof}
 We argue by contradiction. Denote $\theta(\Sigma,\gamma,0)$ by $S$. We discuss the following two cases. 
 
 \textbf{Case 1:} $S>0$.

Suppose the consequence is not true, then $R_\gamma>S$. By 
definition, there is a fill-in $(\Omega,g)$ such that $R_g\geq
S/2$ and $H_g= 0$. Next we construct a metric $g_1$ on $\Omega$ that satisfies $R_{g_1}>S/4$, $g_1|_{\partial\Omega}=\gamma$ and $H_{g_1}>0$ through a conformal deformation. For small $\varepsilon>0$, consider the following equation
\begin{equation*}
\left\{
\begin{aligned}
\Delta_{g}u&=\varepsilon\quad\text{in}\ \Omega,\\
u&=1\quad \text{on}\ \partial\Omega.
\end{aligned}
\right.
\end{equation*}
Obviously, above equation has a unique smooth solution $u$. If $\varepsilon$ is sufficiently small, $u$ is positive. 
Let $g_1=u^\frac{4}{n-2}g$. Then
\begin{align*}
R_{g_1}&=u^{-\frac{n+2}{n-2}}\left(R_{g_1}u-c_n\Delta_{g_1}u\right)\\
&\geq \frac{S}{2}u^{-\frac{4}{n-2}}-c_n\varepsilon u^{-\frac{n+2}{n-2}}.
\end{align*}
Since we have the estimate $|u-1|\leq C\varepsilon$ for some constant $C$, we can find a sufficiently small $\varepsilon$ such that $R_{g_1}>S/4$.   
By the maximum principle, $\frac{\partial
u}{\partial\nu}|_{\partial\Omega}>0$, where $\nu$ is the outward unit normal with respect to $g$. It then follows that
\begin{align*}
H_{g_1}=\frac{c_n}{2}\frac{\partial
u}{\partial\nu}>0.
\end{align*} Thus we obtain a metric $g_1$ on $\Omega$ that satisfies $R_{g_1}>S/4$, $g_1|_{\partial\Omega}=\gamma$ and $H_{g_1}>0$. 

Using Lemma \ref{BMNgluing}, by replacing a neighborhood of $\partial\Omega$ in $(\Omega,g_1)$ with
the ``cylinder" $(\Sigma\times[-\delta,0],\mathrm dt^2+\gamma)$, we can get a
metric $g_2$ on $\Omega$ such that $R_{g_2}>S/8$ and $g_2=\mathrm
dt^2+\gamma$ around a neighborhood of $\partial\Omega$.
Now, we are able to glue the infinitely long cylinder
$(\Sigma\times[0,\infty),\mathrm dt^2+\gamma)$ to $(\Omega,g_2)$ by identifying the slice
$\Sigma\times\{0\}$ with the boundary $\partial\Omega$.
Let $(\tilde\Omega,\tilde g)$ denote the obtained new manifold and
$\tilde\Omega_r$ denote the compact subregion of $\tilde\Omega$ enclosed by
$\Sigma_r:=\Sigma\times\{r\}$.

Take $\phi:[0,+\infty)\to\mathbb R$ to be a smooth function that satisfies
\begin{equation*}
\left\{
\begin{aligned}
&\phi(x)=0\qquad\ \ \   x\in [0,1],\\
&0\leq\phi(x)\leq 1\quad\, x\in [1,2],\\
&\phi(x)=1\qquad\ \ \  x\in [2,+\infty).
\end{aligned}
\right.
\end{equation*}
We also assume $|\phi'|+|\phi''|\leq C$ for some constant $C$. For any $r>0$ and $0<\alpha<1$, define $\tilde u_{r,\alpha}:\tilde\Omega\rightarrow\mathbb{R}$ by
\begin{equation}\label{transition}
\tilde u_{r,\alpha}(x)=\left\{
\begin{aligned}
&\alpha+(1-\alpha)\phi\left(\frac{s}{r}\right) \quad\ \mbox{for}\,\,x\in\Sigma_s,\\
&\alpha \quad\qquad\qquad\qquad\quad\ \ \mbox{elsewhere}.
\end{aligned}
\right.
\end{equation}
It is not hard to see that $\tilde u_{r,\alpha}$ satisfies
\begin{equation*}
\left\{
\begin{aligned}
&\tilde u_{r,\alpha}\equiv \alpha\quad\qquad \text{in}\ \, \tilde\Omega_r, \\
&\alpha\leq \tilde u_{r,\alpha}\leq 1\quad\,\text{in}\ \,\tilde\Omega_{2r}\setminus\tilde\Omega_r,\\
&\tilde u_{r,\alpha}\equiv 1\qquad\quad\, \text{outside}\ \, \tilde\Omega_{2r}, 
\end{aligned}
\right.
\end{equation*}
and
\begin{equation*}
\left|\nabla_{\tilde g}^2 \tilde
u_{r,\alpha}\right|\leq C\left(\alpha^{-\frac{4}{n-2}}r^{-2}+\alpha^{-\frac{n}{n-2}}r^{-1}\right).
\end{equation*}

Make the conformal deformation $\tilde g_1=\tilde
u_{r,\alpha}^{\frac{4}{n-2}}\tilde g$.  A straightforward calculation gives
\begin{equation*}
R_{\tilde g_1}=\left\{
\begin{aligned}
&\alpha^{-\frac{4}{n-2}}R_{\tilde
g}\geq\frac{1}{8}\alpha^{-\frac{4}{n-2}}S \qquad\qquad\qquad\qquad\qquad\qquad\qquad\   \mbox{in}\ \, \tilde\Omega_r,\\
&\tilde u_{r,\alpha}^{-\frac{n+2}{n-2}}\left(R_{\tilde g}\tilde
u_{r,\alpha}-c_n\Delta_{\tilde g}\tilde u_{r,\alpha}\right)\geq R_{\gamma}-c_n\alpha^{-\frac{n+2}{n-2}}\left|\Delta_{\tilde g}\tilde u_{r,\alpha}\right|\ \ \,\ 
\text{in}\ \,\tilde\Omega_{2r}\setminus\tilde\Omega_r,\\
&R_{\gamma}\qquad\qquad\qquad\qquad\qquad\qquad\qquad\qquad\qquad\qquad\qquad\   \mbox{outside}\  \,\tilde\Omega_{2r}.
\end{aligned}
\right.
\end{equation*}

First taking $\alpha$ small enough to guarantee that
$$
\frac{1}{8}\alpha^{-\frac{4}{n-2}}>1,$$
then taking $r$ large enough to ensure that
$$c_n\alpha^{-\frac{n+2}{n-2}}\left|\Delta_{\tilde g}\tilde u_{r,\alpha}\right|<\min
R_\gamma-S,$$
we obtain $\min R_{\tilde g_1}>S$. As a result, $(\tilde\Omega_{3r},
\tilde g_1)$ gives a fill-in of $(\Sigma,\gamma,0)$ with scalar curvature strictly
greater than $S$, which leads to a contradiction.

 \textbf{Case 2:} $S<0$. 

The proof is very similar to the proof for the case $S>0$. 
Suppose the consequence is not true, then $R_\gamma\geq S+\delta$ for some $\delta>0$. By
definition, there exists a fill-in $(\Omega,g)$ such
that $R_g\geq 2S$ and $H_g= 0$. After a conformal deformation similar to that in Case 1, 
we can find a metric $g_1$ on $\Omega$ that satisfies $R_{g_1}\geq 4S$, $g_1|_{\partial\Omega}=\gamma$ and $H_{g_1}>0$. 

Using Lemma \ref{BMNgluing}, by replacing a neighborhood of $\Sigma$ in $(\Omega,g_1)$ with
the ``cylinder" $(\Sigma\times[-\delta,0],\mathrm dt^2+\gamma)$, we can get a
metric $g_2$ on $\Omega$ such that $R_{g_2}>8S$ and $g_2=\mathrm dt^2+\gamma$ around
a neighborhood of $\partial\Omega$.
Now, we are able to glue the infinitely long cylinder
$(\Sigma\times[0,\infty),\mathrm dt^2+\gamma)$ to $(\Omega,g_2)$ by identifying the slice
$\Sigma\times\{0\}$ with the boundary $\partial\Omega$.
Let $(\tilde\Omega,\tilde g)$ denote the obtained new manifold and
$\tilde\Omega_r$ denote the compact subregion of $\tilde\Omega$ enclosed by
$\Sigma_r:=\Sigma\times\{r\}$.

For
any $r>0$ and $\alpha>1$, define $\tilde u_{r,\alpha}:\tilde\Omega\rightarrow\mathbb{R}$ as \eqref{transition}. It is not hard to see that $\tilde u_{r,\alpha}$ satisfies
\begin{equation*}
\left\{
\begin{aligned}
&\tilde u_{r,\alpha}\equiv \alpha\quad\qquad \text{in}\ \, \tilde\Omega_r, \\
&1\leq \tilde u_{r,\alpha}\leq \alpha\quad\,\text{in}\ \,\tilde\Omega_{2r}\setminus\tilde\Omega_r,\\
&\tilde u_{r,\alpha}\equiv 1\qquad\quad\, \text{outside}\ \, \tilde\Omega_{2r}.
\end{aligned}
\right.
\end{equation*}

Make the conformal deformation $\tilde g_1=\tilde
u_{r,\alpha}^{\frac{4}{n-2}}\tilde g$, a straightforward calculation gives
\begin{equation*}
R_{\tilde g_1}=\left\{
\begin{aligned}
&\alpha^{-\frac{4}{n-2}}R_{\tilde
g}\geq 8\alpha^{-\frac{4}{n-2}}S \qquad\qquad\qquad\qquad\qquad\qquad\quad\ \ \,  \mbox{in}\ \, \tilde\Omega_r,\\
&\tilde u_{r,\alpha}^{-\frac{n+2}{n-2}}\left(R_{\tilde g}\tilde
u_{r,\alpha}-c_n\Delta_{\tilde g}\tilde u_{r,\alpha}\right)\geq S+\frac{\delta}{\alpha^4}-c_n\left|\Delta_{\tilde g}\tilde u_{r,\alpha}\right|\ \ \,\ 
\text{in}\ \,\tilde\Omega_{2r}\setminus\tilde\Omega_r,\\
&R_{\gamma}\qquad\qquad\qquad\qquad\qquad\qquad\qquad\qquad\qquad\qquad\quad\ \  \mbox{outside}\  \,\tilde\Omega_{2r}.
\end{aligned}
\right.
\end{equation*}

First taking
$\alpha$ large enough to guarantee that
$$
8\alpha^{-\frac{4}{n-2}}<1,$$
then taking $r$ large enough to ensure that
$$c_n\left|\Delta_{\tilde g}\tilde u_{r,\alpha}\right|\leq\frac{\delta}{2\alpha^4},$$
we obtain $ R_{\tilde g_1}>S$. As a result, $(\tilde\Omega_{3r},\tilde
g_1)$ gives a fill-in of $(\Sigma,\gamma,0)$ with scalar curvature strictly
greater than $S$, which leads to a contradiction.
\end{proof}

\begin{pro}\label{achieved}
If $\theta(\Sigma,\gamma,0)$ can be realized by an extremal fill-in, then
$\theta(\Sigma,\gamma,0)\geq\min R_\gamma$.
\end{pro}

\begin{proof}
By Proposition \ref{neqRbdy}, we only need to consider the case $\theta(\Sigma,\gamma,0)=0$.
We take a contradiction argument.  If the proposition is not true, we may
assume $\min R_\gamma>0$. Suppose
$X:(\Sigma,\gamma)\rightarrow(\Omega,g)$ is an extremal fill-in that
realizes $\theta(\Sigma,\gamma,0)$. By definition, $R_g\geq 
0$ (In fact, by Theorem \ref{staticity}, $R_g\equiv 0$, but we only need $R_g\geq 0$ here). 

Let $l_1$ be an arbitrary positive constant, $l_2=2l_1$, and $l_3$ be a
large constant to be determined. Equip $\Sigma\times[0,l_3]$ with the product
metric $\bar g=\mathrm dt^2+\gamma.$
Glue $(\Sigma\times[0,l_3],\bar g)$ to $(\Omega,g)$ by identifying the slice
$\Sigma\times\{0\}$ with the boundary $\partial\Omega$. It is
obvious that $g=\bar g$ and $H_g=H_{\bar g}=0$ on $\Sigma$. Denote
$\Omega\cup(\Sigma\times[0,l_i])$ by $\Omega_i$ $(i=1,2,3)$. We may modify the differential
structure on $\Omega_3$ so that $(g,\bar g)$ becomes a continuous metric across $\Sigma$. For any
sufficiently small positive constant $\delta$, after carrying out Miao's
mollifying procedure for the metric pair $(g,\bar g)$ (see
\cite{Miao1}), we can get a smooth metric $g_\delta$ on $\Omega_3$
that satisfies:

$\bullet$ $g_\delta=\bar g$ in $\Sigma\times [\delta,l_3]$;

$\bullet$ $g_\delta=g$ in $\Omega\setminus\{\Sigma\times [-\delta,0]\}$;

$\bullet$ $g_\delta$ is uniformly close to $(g,\bar g$) in the
$C^\alpha$-sense for any $0<\alpha<1$;

$\bullet$ $R_{g_\delta}\geq -C$ for some positive $C$ depending only on
$(g,\bar g)$.

Let $R^-_{g_\delta}=\max\{-R_{g_\delta},0\}$. For $u\in
W^{1,2}_0(\Omega_2,g_\delta)$, consider the following functional
\begin{equation*}
I(u)=\int_{\Omega_2}|\nabla_{g_\delta} u|^2-c^{-1}_nR^-_{g_\delta} u^2.
\end{equation*}
By the Sobolev inequality and Minkowski inequality, we have
\begin{align*}
I(u)\geq C_S(\Omega_2,g_\delta)\left\|u\right\|_{L^{\frac{2n}{n-2}}(\Omega_2,g_\delta
)}-c^{-1}_n\left\|u\right\|_{L^{\frac{2n}{n-2}}(\Omega_2,g_\delta
)}\left\|R^-_{g_\delta}\right\|_{L^{\frac{n}{2}}(\Omega_2,g_\delta )},
\end{align*}
where $C_S(\Omega_2,g_\delta)$ is the Sobolev constant of $(\Omega_2,g_\delta)$.

Since $g_\delta$ is uniformly close to $(g,\bar g)$ on $\Omega_2$, for
sufficiently small $\delta$, we have $C_S(\Omega_2,g_\delta)\geq
C_S(\Omega_2,(g,\bar g))/2>0$. On the otherhand,
\begin{equation*}
\left\|R^-_{g_\delta}\right\|_{L^{\frac{n}{2}}(\Omega_2,g_\delta)}=O\big(\delta^{\frac{2}{n}}\big)\quad\text{as}\ \,\delta\rightarrow 0.
\end{equation*}
So for sufficiently small $\delta$ and $u\not\equiv 0$, $I(u)>0$. By the Fredholm alternative
theorem, for any $\varepsilon>0$, the following equation
\begin{equation}\label{makepositive}
\left\{
\begin{aligned}
\Delta_{g_\delta}u+c^{-1}_n R^-_{g_\delta} u&=-\varepsilon\quad\text{in}\ \,
\Omega_2,\\
u&=1\quad\ \ \,\text{on}\ \partial\Omega_2.
\end{aligned}
\right.
\end{equation}
admits a unique solution.
By the Schauder theory, $u\in C^{2,\alpha}(\Omega_2)$. And
$\|u-1\|_{C^{2,\alpha}(\Sigma\times[l_1,l_2])}$ can be arbitrarily small as
$\delta,\varepsilon\rightarrow 0$. Hence, with small $\delta$, $\epsilon$, we can choose a large $l_3$ and get
a positive $C^{2,\alpha}$ function $\tilde u$ on $\overline\Omega_3$ such that $\tilde
u=u$ in $\Omega_1$, $\tilde u\equiv 1$ in a neighborhood of
$\Sigma\times\{l_3\}$ and $$\|\tilde
u-1\|_{C^{2,\alpha}(\Sigma\times[l_1,l_3])}\leq
C\|u-1\|_{C^{2,\alpha}(\Sigma\times[l_1,l_2])},$$
where $C$ is a constant independent of $\delta$ and $\varepsilon$. Set $\tilde
g=\tilde u^{\frac{4}{n-2}}g_\delta$. Then
 \begin{align*}
R_{\tilde g}&=\tilde u^{-\frac{n+2}{n-2}}\left(R_{g_\delta}\tilde
u-c_n\Delta_{g_\delta}\tilde u\right).
\end{align*}
In $\Omega_1$, $\tilde u=u$. By \eqref{makepositive}, we have
 \begin{align*}
R_{\tilde g}&\geq c_nu^{-\frac{n+2}{n-2}}\varepsilon>0.
\end{align*}
In $\Sigma\times[l_1,l_3]$, $g_\delta=\bar g$. It follows that
 \begin{align*}
R_{\tilde g}&\geq\tilde u^{-\frac{4}{n-2}}\min R_\gamma-c_n\tilde
u^{-\frac{n+2}{n-2}}\Delta_{g_\delta}\tilde u.
\end{align*}
So if $\delta$ and $\varepsilon$ are sufficiently small, we have $R_{\tilde g}>0$ 
in $\overline\Omega_3$. Thus $(\Omega_3,\tilde g)$ gives a fill-in of $(\Sigma,\gamma,0)$ with PSC. This contradicts our assumption $\theta(\Sigma,\gamma,0)=0$.
\end{proof}
Having above preparations, we can prove Theorem \ref{nonnegative} in a few words. 
\begin{proof}[Proof of Theorem \ref{nonnegative}]
If $\theta(\Sigma,\gamma,0)\neq 0$, then by Proposition \ref{neqRbdy}, $\theta(\Sigma,\gamma,0)\geq\min R_\gamma>0$. If $\theta(\Sigma,\gamma,0)=0$ and it can be attained, then by Proposition \ref{achieved}, $\theta(\Sigma,\gamma,0)\geq\min R_\gamma>0$. This contradicts the assumption $\theta(\Sigma,\gamma,0)=0$. Consequently, either $\theta(\Sigma,\gamma,0)\geq\min R_\gamma$ or $\theta(\Sigma,\gamma,0)=0$ and it can not be attained.
\end{proof}

We see $R_\gamma>0$ implies $\sigma(\Sigma^{n-1},\gamma,0)\geq 0$. How about the case $R_\gamma$ changes sign but $R_\gamma$ is ``positive on average" in some sense? Inspired by \cite{FS,MS}, we consider the functional $J$ on $W^{1,2}(\Sigma^{n-1})$ defined by
\begin{equation*}
J(f)=\int_{\Sigma}|\nabla_{\gamma} f|^2+\frac{1}{2}R_\gamma f^2\,\mathrm d\mu_\gamma.
\end{equation*}
Let $\lambda_1$ be its first eigenvalue and $f_1$ be the corresponding
eigenfunction.
Then $f_1>0$ and satisfies
\begin{equation*}
-\Delta_\gamma f_1+\frac{R_\gamma}{2}f_1=\lambda_1f_1.
\end{equation*}
The metric $\bar g_1=f_1^2dt^2+\gamma$ has constant scalar curvature $R_{\bar g_1}\equiv 2\lambda_1$. Obviously, $2\lambda_1(\gamma)\geq\min R_\gamma$. 
Using the cylinder $(\Sigma^{n-1}\times I,\bar g_1)$ in suitable gluing constructions, we can get the stronger version of Theorem \ref{nonnegative}, namely
\begin{thm}\label{strongnonnegative}
For  $(\Sigma^{n-1},\gamma)$ with $\lambda_1(\gamma)>0$, either
\begin{enumerate}
\item[(1)] $\theta(\Sigma^{n-1},\gamma,0)\geq 2\lambda_1(\gamma)$, or
\item[(2)] $\theta(\Sigma^{n-1},\gamma,0)=0$ and it can not be attained.
\end{enumerate}
\end{thm}

\subsection{Existence of fill-in with PSC metrics}\ 
\vskip 0.1cm
In this subsection, by virtue of the Schwarzschild neck, we extend the results of minimal fill-ins. We first prove Theorem \ref{fillin1} by utilizing the monotonicity of $\theta$-invariant (see Theorem \ref{Thm: weak monotonicity sigma
invariant}), Proposition \ref{neck}, Proposition \ref{neqRbdy} and Proposition \ref{achieved}. 

\begin{proof}[Proof of Theorem \ref{fillin1}]

Since $R_\gamma>0$, from Theorem \ref{nonnegative}, we know that either 
$\theta(\Sigma^{n-1},\gamma, 0)\geq\min R_\gamma$,  
or $\theta(\Sigma^{n-1},\gamma,0)=0$ and it can not be
attained. If $\theta(\Sigma^{n-1},\gamma, 0)\geq \min R_\gamma$, then by the definition
of $\theta$-invariant, for any $\varepsilon>0$, there exists a fill-in $(\Omega_1,g_1)$ of
$(\Sigma^{n-1},\gamma,0)$ with $R_{g_1} >\min R_\gamma-\varepsilon$. Choose $\varepsilon<\min R_\gamma$. According to Proposition \ref {neck}, for any positive constant $h<\max H$, there exists a Schwarzschild neck {$(\Sigma\times
[r_1,r_2],g)$} with 
\begin{equation}\label{scalarlowerbound}
R_g\geq \min R_\gamma-\frac{n-2}{n-1}\max H^2-\varepsilon,
\end{equation}
 whose boundary data is $(\Sigma^{n-1},\gamma,\max H)$ and  
$(\Sigma^{n-1},\mu\gamma,h)$. Here, $\mu=r_1^2\psi^2(r_1)<1$.
Moreover, by the property of Schwarzschild neck, we can extend
$(\Sigma^{n-1}\times [r_1,r_2],g)$ a little to obtain $(\Sigma^{n-1}\times [r'_1,r_2],g)$
($r'_1<r_1$). And the scalar curvature of the extended Schwarzschild neck still satisfies \eqref{scalarlowerbound}. In a similar manner to the gluing construction in the proof of Proposition \ref{Prop: strict inequality 2}, we glue $(\Sigma^{n-1}\times
[r_1,r_2],g)$ to $(\Omega_1, \mu g_1)$ along the slice {$\Sigma^{n-1}\times
\{r_1\}$} to get a fill-in $(\Omega_2,g_2)$ of $(\Sigma^{n-1},\gamma,\max H)$, which satisfies
$R_{g_2}\geq \min R_\gamma-\frac{n-2}{n-1}\max H^{2}-2\varepsilon$ (here we need
to use the fact that $\mu<1$). Since $\varepsilon$ can be arbitrarily small, in fact we have 
$$\theta(\Sigma^{n-1},\gamma,\max H)\geq \min
R_\gamma-\frac{n-2}{n-1}\max H^{2}.$$ By the monotonicity of $\theta$-invariant, 
$$\theta(\Sigma^{n-1},\gamma,H)\geq\theta(\Sigma^{n-1},\gamma,\max H).$$
Combing above two inequalities together, we get the desired inequality. 

Consider the second case: $\theta(\Sigma^{n-1},\gamma,0)=0$ and $\theta(\Sigma^{n-1},\gamma,0)$ can not be
attained. By definition, for any $\varepsilon>0$, there exists a fill-in $(\Omega_3,g_3)$ of $(\Sigma^{n-1},\gamma,0)$ with $R_{g_3}\geq -\varepsilon$. After a similar gluing construction as above, we can obtain a fill-in $(\Omega_4,g_4)$ of
$(\Sigma^{n-1},\gamma,\max H)$ with $R_{g_4}\geq -C\varepsilon$ for a constant $C$ indepentdent of $\varepsilon$. Since $\varepsilon$ can be arbitrarily small, actually we have $\theta(\Sigma^{n-1},\gamma,\max H)\geq 0$. On the other hand, the monotonicity formula tell us that $$\theta(\Sigma^{n-1},\gamma,\max H)\leq\theta(\Sigma^{n-1},\gamma,H)\leq
\theta(\Sigma^{n-1},\gamma,0).$$ Therefore, $\theta(\Sigma^{n-1},\gamma,H)=0$. 
Next we prove $\theta(\Sigma^{n-1},\gamma,H)$ can not be
attained by a contradiction argument. Suppose $(\Omega,g)$ is an extremal fill-in that
realizes $\theta(\Sigma^{n-1},\gamma,H)$. By definition, $R_g\geq
0$. Then we glue  a very long cylinder $\Sigma\times [0,l]$
equipped with the metric $\mathrm dt^2+\gamma$ to $(\Omega,g)$ through Miao's gluing procedure for metrics with corners. Then by a very similar argument as the
proof of Proposition \ref{achieved}, we get $\theta(\Sigma^{n-1},\gamma,0)>0$, 
which contradicts to the assumption $\theta(\Sigma^{n-1},\gamma,H)=0$. This completes the proof. 
\end{proof}

In the sequel, we give the proof of Theorem \ref{Thm: main1} in two steps. In the first step, we prove Theorem
\ref{Thm: main1} for the case that $H$ is an arbitrary negative constant. In the second step, we apply Proposition \ref{neck} to deform $(\Sigma^{n-1},\gamma_0,H)$ to $(\Sigma^{n-1},\mu\gamma_0,-\epsilon)$
for some positive constants $\mu$ and $\epsilon$, and the deformation
provides a manifold with PSC. Using a gluing argument from
Lemma \ref{BMNgluing}, we finally get the desired result.

We are going to show the following:

\begin{pro}\label{Lem: connection}
Let $\gamma_0$ and $\gamma_1$ be two smooth metrics in $\mathcal
M^k_{psc}(\Sigma^{n-1})$ isotopy to each other. Given any $\epsilon_0>0$, there exist positive constants $\mu$ and $\epsilon_1$
 such that $(\Sigma^{n-1},\gamma_0,-\epsilon_0)$ and
$(\Sigma^{n-1},\mu\gamma_1,\epsilon_1)$ can be realized as the boundary data of $(\Sigma^{n-1}\times[0,1],g)$ for some metric $g$ with positive scalar curvature.
\end{pro}
\begin{proof}
By assumption, there is a continuous path $\{\gamma(t)\}_{t\in[0,1]}$ in $\mathcal{M}^k_{psc}(\Sigma^{n-1})$ with $\gamma(0)=\gamma_0$ and $\gamma(1)=\gamma_1$. Without loss of generality, we may assume $\gamma_t\equiv\gamma_0$ around $t=0$ and $\gamma_t\equiv \gamma_1$ around $t=1$.
By Proposition 2.1 and its proof in \cite{CAM}, we may also assume the path is smooth. 
Define a metric $g$ on $\Sigma\times[0,1]$ by
\begin{equation*}
g=A^2(t)\mathrm dt^2+e^{2B(t)}\gamma(t), 
\end{equation*}
with functions $A$ and $B$ on $[0,1]$ to be determined later.
Define another metric $\bar g$ on $\Omega$ by 
\begin{equation*}
\bar g=\mathrm dt^2+e^{2B(t)}\gamma(t).
\end{equation*}
A straightforward calculation gives the mean curvature of the slice
$\Sigma_t=\Sigma\times\{t\}$ with respect to the metric $\bar g$ and the
 $\partial_t$-direction, namely
\begin{equation*}
\bar H_t=(n-1)B'(t)+\frac{1}{2}\tr_{\gamma(t)}\gamma'(t).
\end{equation*}
We may assume that $|\tr_{\gamma(t)}\gamma'(t)|\leq C_1$ for some positive constant $C_1$. Note that $\gamma'(t)\equiv 0$ around $t=0$ and $t=1$. 
 Fixing $c_0\in(0,\epsilon_0)$, we can choose
$B$ such that $\bar H_0\equiv\bar H_1\equiv \epsilon_0$
and $\bar H_t\geq c_0$ for all $t\in[0,1]$. When the function $B$ is chosen, we may assume $|R_{\bar g}|<C_2$ for some positive constant $C_2$. From the calculations in \cite{Shi-Tam1}, we know
\begin{equation*}
R_{g}=2\bar H_tA^{-3}A'+\left(1-A^{-2}\right) R_t+A^{-2}R_{\bar g},
\end{equation*}
where $R_t$ is the scalar curvature of $(\Sigma_t,e^{2B(t)}\gamma_t)$. Take
$A(t)=e^{\Lambda t}$ with $\Lambda$ to be determined. Then
\begin{equation*}
R_g=e^{-2\Lambda t}\left(2\Lambda \bar H_t+(e^{2\Lambda
t}-1)R_t+R_{\bar g}\right).
\end{equation*}
Since $\bar H_t\geq c_0>0$,  $R_t>0$ and $|R_{\bar g}|\leq C_2$, we may choose sufficiently large $\Lambda$ so that $R_g>0$. Now,
it is clear that the mean curvature of $\Sigma_0$ is $$H_{\Sigma_0}=-A(0)^{-1}\bar H_0=-\epsilon_0,$$
and the mean curvature of $\Sigma_1$ is $$H_{\Sigma_1}=A(1)^{-1}\epsilon_0=e^{-\Lambda}\epsilon_0.$$
It is easy to see that $\mu=e^{2B(1)}$ and $\epsilon_1=e^{-\Lambda}\epsilon_0$.
\end{proof}

\begin{pro}\label{Prop: fill-in of connected metric}
Let $\gamma_0$ and $\gamma_1$ be two smooth metrics in $\mathcal
M^k_{psc}(\Sigma^{n-1})$ isotopic to each other. If
$(\Sigma^{n-1},\gamma_1,0)$ admits a fill-in of positive scalar curvature,
then $(\Sigma^{n-1},\gamma_0,-\epsilon_0)$ admits a fill-in of positive
scalar curvature for any $\epsilon_0>0$.
\end{pro}
\begin{proof}
By Proposition \ref{Lem: connection}, for any $\epsilon_0>0$, there exists a smooth metric $g$ on
$\Sigma\times[0,1]$ such that $R_g>0$ and the boundary data is
$(\Sigma,\gamma_0,-\epsilon_0)$ and $(\Sigma,\mu\gamma_1,\epsilon_1)$ for some positive constants $\mu$ and $\epsilon_1$. Denote $\Sigma\times[0,1]$ by $\Omega$. From the construction of $g$, it is clear that $g$ has the form
\begin{equation*}
g=A^2(t)\mathrm dt^2+e^{2B(t)}\gamma_1
\end{equation*}
in the Gaussian coordinate of a small $\delta$-collar neighborhood of $\Sigma\times\{1\}$ (corresponding to $(\Sigma,\mu\gamma_1,\epsilon_1)$) in $\Omega$.
Let $r=\int_{1-t}^1A(s)\,\mathrm d s$.
This collar neighborhood can be represented by
$[0,\delta]\times\Sigma$ and $g$ can be represented by
\begin{equation*}
g=\mathrm dr^2+q(r)\gamma_1,
\end{equation*}
where $q$ is a positive function on $[0,\delta]$ with $q(0)=\mu$. 
Since $(\Sigma,\gamma_1,0)$ admits a fill-in of PSC, $(\Sigma,q(\delta)\gamma_1,0)$ admits a fill-in of PSC, denoted by $(\tilde\Omega,\tilde g)$. 
Note that
$\{\delta\}\times\Sigma$ has negative mean curvature with respect to the $\partial_r$-direction. In a similar manner to the gluing construction in the proof of Proposition \ref{Prop: strict inequality 2}, we can glue the complement of the $\delta$-collar neighborhood in $\Omega$ to $(\tilde\Omega,\tilde g)$ along the slice $\{\delta\}\times\Sigma$ to obtain a new manifold
of PSC, whose boundary data is
$(\Sigma,\gamma_0,-\epsilon_0)$. This gives the desired fill-in. 
\end{proof}

After above preparations, we give the proof of Theorem \ref{Thm: main1}.
\begin{proof}[Proof of Theorem \ref{Thm: main1}]
According to Theorem \ref{Thm: weak monotonicity sigma invariant}, it suffices to
prove the theorem for positive $H$. By Proposition \ref{neck}, we can find a Schwarzschild neck
$(\Omega,g)$ with $R_g>0$, whose boundary data is $(\Sigma,\gamma_0,H)$ and
$(\Sigma,\mu\gamma_0,\epsilon)$ with $\mu,\epsilon>0$. By the property of Schwarzschild neck, we can extend
$(\Omega,g)$ through $(\Sigma,\mu\gamma_0,\epsilon)$ a little bit, and the extended neck still has PSC. 
Since $\gamma_0$ is isotopic to $\gamma_1$ in $\mathcal
M^k_{psc}(\Sigma)$, $\mu\gamma_0$ is isotopic to $\mu\gamma_1$ in $\mathcal
M^k_{psc}(\Sigma)$. It then follows from
Proposition \ref{Prop: fill-in of connected metric} that
$(\Sigma,\mu\gamma_0,-\epsilon/2)$ admits a fill-in $(\tilde\Omega,\tilde
g)$ of PSC. In a similar manner to the gluing construction in the proof of Proposition \ref{Prop: strict inequality 2}, we can glue $(\Omega,g)$ to $(\tilde\Omega,\tilde g)$ along $(\Sigma,\mu\gamma_0)$ to obtain a fill-in of $(\Sigma,\mu\gamma_0,H)$ with PSC.
\end{proof}

In fact, we can strengthen Theorem \ref{Thm: main1} to the following one:
\begin{thm}
Let $\gamma_0$ and $\gamma_1$ be two smooth metrics in $\mathcal
M^k _{psc}(\Sigma^{n-1})$ isotopic to each other. If
$(\Sigma^{n-1},\gamma_1,-H_1)$ admits a fill-in with positive scalar curvature for
some positive function $H_1$ satisfying
\begin{equation*}
H_1<\left(\frac{n-1}{n-2}\min R_{\gamma_1}\right)^{\frac{1}{2}},
\end{equation*}
then for any function $H$ with
\begin{equation*}
H<\left(\frac{n-1}{n-2}\min R_{\gamma_0}\right)^{\frac{1}{2}},
\end{equation*}
$(\Sigma^{n-1},\gamma_0,H)$ admits a fill-in with positive scalar curvature.
\end{thm}
From Theorem \ref{Thm: main1}, it suffices to show the following:
\begin{lm}
Let $(\Sigma^{n-1},\gamma,-H)$ be a triple of Bartnik data. Suppose $H$ is a positive constant and $R_\gamma>\frac{n-2}{n-1}H^2$.
If $(\Sigma^{n-1},\gamma,-H)$ admits a fill-in of positive scalar curvature, then
$(\Sigma^{n-1},\gamma,0)$ admits a fill-in of positive scalar curvature.
\end{lm}
\begin{proof}
Suppose $(\Omega_1,g_1)$ is a fill-in of $(\Sigma^{n-1},\gamma,-H)$ with
PSC. Since $R_\gamma>\frac{n-2}{n-1}H^2$, there is a constant $\delta>0$ such that
$$
 R_\gamma>\frac{n-2}{n-1}(H+\delta)^2.
$$
By Proposition \ref{neck}, we can find a Schwarzschild neck $(\Omega_2,g_2)$
with PSC, whose boundary data is
$(\Sigma^{n-1},\gamma,H+\delta)$ and $(\Sigma^{n-1},\mu\gamma,0)$  for a certain positive constant $\mu$. Note that we can extend
$(\Omega_2,g_2)$ through $(\Sigma^{n-1},\gamma,H+\delta)$ a little bit. In a similar manner to the gluing construction in the proof of Proposition \ref{Prop: strict inequality 2}, we may glue $(\Omega_2,g_2)$ to $(\Omega_1,g_1)$ along $(\Sigma^{n-1},\gamma)$, to obtain a fill-in of $(\Sigma^{n-1},\mu\gamma,0)$ with PSC. The desired fill-in then follows by a rescaling.
\end{proof}

{\it Acknowledgements.} We would like to thank M. Gromov,  J. L. Jauregui, P. Miao and C. Sormani for their interest in this work. We are thankful to C. Mantoulidis for his suggestion to make the proof of Theorem \ref{nonexistfillin1} clearer and interest in this work.

\end{document}